\documentclass{article}%
\usepackage{amsfonts}
\usepackage{amsmath}
\usepackage{amssymb}
\usepackage{amsxtra}
\usepackage{graphicx}
\usepackage{geometry}
\usepackage{color}
\usepackage{colortbl}%
\providecommand{\U}[1]{\protect\rule{.1in}{.1in}}
\newtheorem{theorem}{Theorem}

\newtheorem{corollary}[theorem]{Corollary}

\newtheorem{lemma}[theorem]{Lemma}

\newtheorem{remark}[theorem]{Remark}

\numberwithin{equation}{section}
\begin{document}

\title{Water-waves modes trapped in a canal by a body \\with the rough surface}
\author{G.Cardone\\University of Sannio - Department of Engineering\\Corso Garibaldi, 107 - 82100 Benevento, Italy\\email: giuseppe.cardone@unisannio.it
\and T.Durante\\University of Salerno\\Department of Information Engineering and Applied Mathematics\\Via Ponte don Melillo, Fisciano (SA), Italy\\email: durante@diima.unisa.it
\and S.A.Nazarov\\Institute of Mechanical Engineering Problems\\V.O., Bolshoi pr., 61, 199178, St. Petersburg, Russia.\\email: srgnazarov@yahoo.co.uk}
\maketitle

\begin{abstract}
The problem about a body in a three dimensional infinite channel is considered
in the framework of the theory of linear water-waves. The body has a rough
surface characterized by a small parameter $\varepsilon>0$ while the distance
of the body to the water surface is also of order $\varepsilon$. Under a
certain symmetry assumption, the accumulation effect for trapped mode
frequencies is established, namely, it is proved that, for any given $d>0$ and
integer $N>0$, there exists $\varepsilon(d,N)>0$ such that the problem has at
least $N$ eigenvalues in the interval $(0,d)$ of the continuous spectrum in
the case $\varepsilon\in\left(  0,\varepsilon(d,N)\right)  $. The
corresponding eigenfunctions decay exponentially at infinity, have finite
energy, and imply trapped modes.



\end{abstract}

\section{Introduction}

\subsection{Statement of the problem.}

\begin{figure}
[ptb]
\begin{center}
\includegraphics[
height=1.3258in,
width=2.9205in
]%
{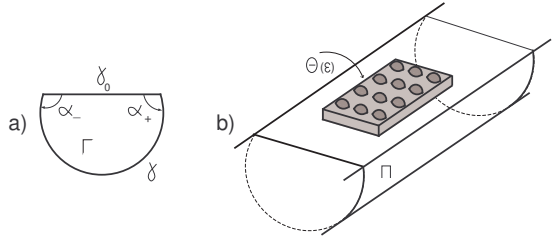}%
\caption{The transverse cross-section $\Gamma$ and the infinite cylindrical
canal with submerged body $\Theta\left(  \varepsilon\right)  $}%
\label{f1}%
\end{center}
\end{figure}

Let $\Gamma$ be a domain on the plane $\mathbb{R}^{2}\ni x^{\prime}=\left(
x_{2,}x_{3}\right)  $ bounded by the line interval $\mathbf{\gamma}%
_{0}=\left\{  x^{\prime}:\left\vert x_{2}\right\vert <l,\text{\ }%
x_{3}=0\right\}  $ and the smooth simple curve $\mathbf{\gamma}$ inside the
lower half-plane $\mathbb{R}_{-}^{2}=\left\{  x^{\prime}:x_{3}<0\right\}  $
which meets $\mathbf{\gamma}_{0}$\textbf{\ }at the points $x^{\prime}=\left(
\pm l,0\right)  $ with the angles $\alpha_{\pm}\in\left(  0,\pi\right)  $ (see
Fig. \ref{f1}, a).%

\begin{figure}
[ptb]
\begin{center}
\includegraphics[
height=0.9928in,
width=1.6371in
]%
{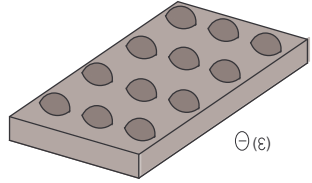}%
\caption{The body $\mathbf{\Theta}\left(  \varepsilon\right)  $ with the upper
rugged surface.}%
\label{f2}%
\end{center}
\end{figure}

\begin{figure}
[ptb]
\begin{center}
\includegraphics[
height=1.107in,
width=3.0865in
]%
{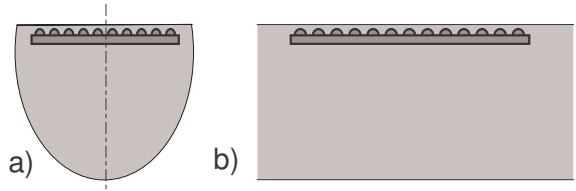}%
\caption{The tranverse (a) and longitudinal (b) cross-sections of the body
$\mathbf{\Theta}\left(  \varepsilon\right)  $ in the canal.}%
\label{f3}%
\end{center}
\end{figure}

\begin{figure}
[ptb]
\begin{center}
\includegraphics[
height=1.2851in,
width=3.0675in
]%
{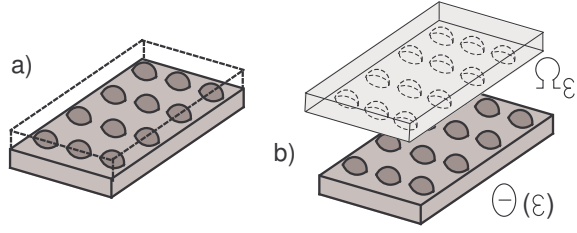}%
\caption{The body $\mathbf{\Theta}\left(  \varepsilon\right)  $ and the
rectangular plate-shaped layer $\Omega_{\varepsilon}$ of water, jointed (a)
and separated (b).}%
\label{f4}%
\end{center}
\end{figure}

The three-dimensional canal $\Pi=\mathbb{R}\mathbf{\times}\Gamma\ni x=\left(
x_{1},x^{\prime}\right)  $ with the horizontal plain surface $\Lambda
=\mathbb{R}\mathbf{\times\gamma}_{0}$ contains the finite body\textbf{\ }%
$\mathbf{\Theta}\left(  \varepsilon\right)  $ (see Fig. \ref{f1}, b). The
shape of the body depends on the small parameter $\varepsilon>0$ so that its
upper surface is rough with periodic fine knobs and/or caverns of size
$\varepsilon$ (see Fig. \ref{f2} with the three-dimensional image and Fig.
\ref{f3} with the two-dimensional cross-sections). The body is submerged in
the superficial region and the mean distance from $\Lambda$ to the upper
surface of $\mathbf{\Theta}\left(  \varepsilon\right)  $ is of order
$\varepsilon$ as well. There is no geometrical restriction on the bottom of
$\mathbf{\Theta}\left(  \varepsilon\right)  $ but the upper rough horizontal
part of the boundary $\partial\mathbf{\Theta}\left(  \varepsilon\right)  $
restricts from below a finite thin rectangular plate-shaped part
$\Omega_{\varepsilon}$ of the near-surface water layer (see Fig. \ref{f4})%
\begin{equation}
\Omega_{\varepsilon}\subset\Pi\left(  \varepsilon\right)  =\Pi\backslash
\overline{\mathbf{\Theta}\left(  \varepsilon\right)  } \label{1.1}%
\end{equation}
In other words, the upper straight base $\omega_{+}$ of the plate
$\Omega_{\varepsilon}$ belongs to the horizontal surface $\Lambda$ of water
while the lower base $\omega_{-}\left(  \varepsilon\right)  $ of a fine
periodic structure reposes upon the boundary $\partial\mathbf{\Theta}\left(
\varepsilon\right)  $ of the body.%

%

%

Although further results are valid for Lipschitz surface $\omega_{-}\left(
\varepsilon\right)  $ (see Section \ref{sect4}), we assume in the presentation
that this surface is smooth enough though. The assumption crucially simplifies
rather cumbersome calculations in Sections \ref{sect2.4} and \ref{sect2.5}.

To describe the periodic structure of the plate $\Omega_{\varepsilon}$ more
precisely, we introduce the periodicity cell $\Sigma$ such that%
\[
\sigma\times\left(  -h,0\right)  \subset\Sigma\subset\sigma\times\left(
-H,0\right)  ,
\]
where $\sigma=\left\{  y=\left(  y_{1},y_{2}\right)  \in\mathbb{R}%
^{2}:\left\vert y_{i}\right\vert <a_{i}/2,\text{ }i=1,2\right\}  $ is a
rectangle ($a_{i}>0)$ and $0<h\leq H.$ We introduce another rectangle
\begin{equation}
\omega=\left\{  y:\left\vert y_{i}\right\vert <A_{i}/2,\text{ }i=1,2\right\}
\label{1.2}%
\end{equation}
and assume that the sizes are in the relation
\[
A_{i}=\varepsilon a_{i}N_{i},\text{ \ }i=1,2,
\]
where $N_{1}$ and $N_{2}$ are large positive integers. We then set%
\begin{equation}
\Sigma_{\varepsilon}^{\nu}=\left\{  x=\left(  y,z\right)  :\left(
\varepsilon^{-1}y_{1}-\nu_{1}a_{1},\varepsilon^{-1}y_{2}-\nu_{2}%
a_{2},\varepsilon^{-1}z\right)  \in\Sigma\right\}  , \label{1.3}%
\end{equation}%
\begin{equation}
\overline{\Omega}_{\varepsilon}=\underset{\nu:\left\vert \nu_{i}\right\vert
\leq N_{i}}{%
{\displaystyle\bigcup}
}\overline{\Sigma_{\varepsilon}^{\nu}}, \label{1.4}%
\end{equation}
where $\nu=\left(  \nu_{1},\nu_{2}\right)  \in\mathbb{Z}^{2}$ is a multi-index
and $\mathbb{Z}\mathbf{=}\left\{  0,\pm1,....\right\}  .$ The domain
$\Omega_{\varepsilon}$, i.e., the interior of the closed set (\ref{1.3}) is
but a thin plate composed from the large number $N_{1}\times N_{2}$ of small
periodicity cells (\ref{1.3}). We do not exclude the case $h=H$ when $\Sigma$
and $\Pi$ imply parallelepipedes of sizes $a_{1}\times a_{2}\times h$ and
$A_{1}\times A_{2}\times\varepsilon H,$ respectively (cf. \cite{NaMatSb}).

Notice that it is convenient to use different notation for the same Cartesian
coordinate system $\ x=\left(  x_{1},x_{2},x_{3}\right)  ,$ namely, $\left(
x_{1,}x^{\prime}\right)  $ in the canal $\Pi\left(  \varepsilon\right)  $ and
$\left(  y,z\right)  $ in the plate $\Omega_{\varepsilon}$ while $x^{\prime
}=\left(  x_{2},x_{3}\right)  $ are coordinates on the cross-section $\Gamma$
of $\Pi$ and $y=\left(  y_{1},y_{2}\right)  =\left(  x_{1},x_{2}\right)  $ are
coordinates on the upper base $\omega_{+}=\left\{  x:y\in\omega,\text{
}z=0\right\}  $ of $\Omega_{\varepsilon}.$

In the canal $\Pi$ with the submerged body $\mathbf{\Theta}\left(
\varepsilon\right)  ,$ we consider the spectral problem of the linearized
water-wave theory
\begin{equation}
-\Delta_{x}\Phi_{\varepsilon}\left(  x\right)  =0,\text{ \ }x\in\Pi\left(
\varepsilon\right)  , \label{1.5}%
\end{equation}%
\begin{equation}
\partial_{n}\Phi_{\varepsilon}\left(  x\right)  =0,\text{ \ }x\in\pi\left(
\varepsilon\right)  :=\partial\Pi\left(  \varepsilon\right)  \backslash
\overline{\Lambda}, \label{1.6}%
\end{equation}%
\begin{equation}
\partial_{z}\Phi_{\varepsilon}\left(  x\right)  =\lambda_{\varepsilon}%
\Phi\left(  x\right)  ,\text{ \ }x\in\Lambda\label{1.7}%
\end{equation}
(see, e.g., monographs \cite{KMV, Stoker, John} for physical and mathematical
background). Here $\Delta_{x}=\nabla_{x}\cdot\nabla_{x}$ is the Laplace
operator, $\nabla_{x}=\operatorname{grad}$ and $\nabla_{x}\cdot$
$=\ \operatorname{div}$, while $\partial_{n}$ is the derivative along the
outward normal, in particular, $\partial_{n}=\partial_{z}$ on $\Lambda.$
Furthermore, $\Phi_{\varepsilon}$ is the velocity potential and $\lambda
_{\varepsilon}$ the spectral parameter, proportional to square of frequency of
harmonic oscillations in the canal.

In addition to the smoothness assumptions introduced above, the whole boundary
$\partial\Pi\left(  \varepsilon\right)  $ is Lipschitz. Hence, the normal $n$
and the Neumann (\ref{1.6}) and the Steklov (\ref{1.7}) boundary conditions
are defined properly for almost all $x\in\partial\Pi\left(  \varepsilon
\right)  .$ However, the gradient $\nabla_{x}\Phi_{\varepsilon}$ can gain
singularities, e.g., at edges on the boundary and in Section \ref{sect1.3} we
give a precise definition of an operator $\mathcal{L}_{\varepsilon}$ for
problem (\ref{1.5})-(\ref{1.7}) in the Sobolev space $H^{1}\left(  \Pi\left(
\varepsilon\right)  \right)  .$ In this framework, being interested to detect
trapped modes, i.e., solutions with the exponential decay as $x_{1}%
\rightarrow\pm\infty,$ we need not to supply the problem with any radiation
condition at infinity. We again refer to \cite{KMV,Stoker,John}\ for
formulation of these radiation conditions in similar geometrical situations.

\subsection{The trapped modes frequencies}

In this paper we seek for trapped modes, i.e., solutions $\Phi_{\varepsilon
}\in H^{1}(\Pi(\varepsilon))$ of problem (\ref{1.5})-(\ref{1.7}) with the
finite energy and, therefore, the exponential decay at infinity. Such
solutions have been a goal in many investigations (see \cite{N1}-\cite{N9} and
review \cite{ML} for much more extensive list of references). In the sequel we
detect the accumulation effect of trapped mode eigenfrequencies, namely,
assuming the geometrical parameter $\varepsilon>0$ sufficiently small, we find
out any prescribed number $N$ of eigenvalues on the given small interval
$\left(  0,d\right)  ,$ $d>0,$ of the continuous spectrum in problem
(\ref{1.5})-(\ref{1.7}). We make use of the following issues:

\begin{itemize}
\item The artificial Dirichlet boundary conditions on the plane $\left\{
x:x_{1}=0\right\}  .$

\item Asymptotic analysis for eigenvalues of a spectral problem in the thin
finite domain $\Omega_{\varepsilon}.$

\item The operator formulation of the problem in Hilbert space.

\item The max-min principle.
\end{itemize}

Let us outline these issues.

First, the artificial Dirichlet boundary conditions on the plane of
geometrical symmetry permit to create a positive threshold $\lambda\left(
\Gamma\right)  >0$ in the modified spectral problem so that the continuous
spectrum covers the ray $\left[  \lambda\left(  \Gamma\right)  ,+\infty
\right)  $ but leaves the gap $\left(  0,\lambda\left(  \Gamma\right)
\right)  $ for the discrete spectrum. This trick was proposed in \cite{Vas}
for detecting trapped modes in a strip with a symmetric obstacle for the
Helmholts equation with the Neumann boundary conditions.

Second, as a subsidiary problem, we investigate sloshing mode eigenfrequencies
in the artificially constructed thin finite layer $\Omega_{\varepsilon}$ of
water (see formula (\ref{1.1}) and Fig. \ref{f4} where it is demonstrated how
the plate-shaped layer $\Omega_{\varepsilon}$ is cut off and separated by the
body $\Theta\left(  \varepsilon\right)  ).$ In other words, we consider the
auxiliary Steklov spectral problem%
\begin{equation}
-\Delta_{x}u_{\varepsilon}\left(  x\right)  =0,\text{ \ }x\in\Omega
_{\varepsilon}, \label{1.12}%
\end{equation}%
\begin{equation}
\partial_{n}u_{\varepsilon}\left(  x\right)  =0,\text{ \ \ \ }x\in\omega
_{-}\left(  \varepsilon\right)  ,\text{ \ }\partial_{z}u_{\varepsilon}\left(
x\right)  =\alpha_{\varepsilon}u_{\varepsilon}\left(  x\right)  ,\text{
\ \ \ }x\in\omega_{+}, \label{1.13}%
\end{equation}%
\begin{equation}
u_{\varepsilon}\left(  x\right)  =0,\text{ \ }x\in\Upsilon_{\varepsilon
}=\partial\Omega_{\varepsilon}\backslash\left(  \overline{\omega^{-}\left(
\varepsilon\right)  \cup\omega^{+}}\right)  . \label{1.14}%
\end{equation}
The asymptotic analysis of (\ref{1.12})-(\ref{1.14}) is rather standard (cf.
\cite{BLP, SP, Nabook} and others). However, a new effect is observed in
Theorem \ref{Theorem 2.11}: each entry of the monotone unbounded eigenvalue
sequence in problem (\ref{1.12})-(\ref{1.14})
\begin{equation}
0<\alpha_{\varepsilon}^{\left(  1\right)  }<\alpha_{\varepsilon}^{\left(
2\right)  }\leq...\leq\alpha_{\varepsilon}^{\left(  N\right)  }\leq
....\rightarrow+\infty\label{1.88}%
\end{equation}
becomes infinitesimal when $\varepsilon\rightarrow0^{+}.$ Namely, the
eigenvalues $\alpha_{\varepsilon}^{\left(  1\right)  },\alpha_{\varepsilon
}^{\left(  2\right)  },...,\alpha_{\varepsilon}^{\left(  N\right)  }$ belong
to the interval $\left(  0,d\right)  \subset\left(  0,\lambda\left(
\Gamma\right)  \right)  $ if $\varepsilon<\varepsilon\left(  d,N\right)  ,$
with a certain $\varepsilon\left(  d,N\right)  >0.$

It suffices to prove that the point spectrum of problem (\ref{1.5}%
)-(\ref{1.7}) in the interval $\left(  0,d\right)  $ contain at least $N$
eigenvalues. This task is fulfilled by applying the max-min principle (see,
e.g., \cite[Theorem 10.2.2]{BiSo}) to the operator formulation \cite{NaMatSb}
of the problem (respectively the fourth and third issues in the above list).
We emphasize that the lateral side $\Upsilon_{\varepsilon}$ of the plate
$\Omega_{\varepsilon}$ is supplied with the Dirichlet conditions (that is why
we call (\ref{1.12})-(\ref{1.14}) the Steklov spectral problem while the
complete analogy with sloshing modes is dubious). We again use the geometrical
symmetry and reduce the problem (\ref{1.12})-(\ref{1.14}) onto the subdomain
$\Omega_{\varepsilon}^{+}=\left\{  x\in\Omega_{\varepsilon}:x_{2}>0\right\}
.$ Imposing the Dirichlet condition on the artificial boundary $\left\{
x\in\Omega_{\varepsilon}:x_{2}=0\right\}  ,$ we keep the concentration
property for eigenvalues $\alpha_{\varepsilon}^{\left(  p\right)  +}$ of the
Steklov problem in $\Omega_{\varepsilon}^{+}.$ The Dirichlet conditions and
the inclusion $\omega^{-}\left(  \varepsilon\right)  \subset\partial
\mathbf{\Theta}\left(  \varepsilon\right)  $ permit for the extension of the
corresponding eigenfunctions $u_{\varepsilon}^{\left(  p\right)  +}$ by zero
from $\Omega_{\varepsilon}^{+}$ onto the set%
\begin{equation}
\Pi^{+}\left(  \varepsilon\right)  =\left\{  x\in\Pi\left(  \varepsilon
\right)  :x_{2}>0\right\}  . \label{1.10}%
\end{equation}
These extended eigenfunctions are taken as trial functions in the max-min
principle which ensure that, for any $\alpha_{\varepsilon}^{\left(  p\right)
+}\in\left(  0,\lambda\left(  \Gamma\right)  \right)  ,$ the point spectrum of
the problem in $\Pi^{+}\left(  \varepsilon\right)  $ contains an eigenvalue
$\lambda_{\varepsilon}^{\left(  p\right)  +}\in\left(  0,\alpha_{\varepsilon
}^{\left(  p\right)  +}\right]  \subset\left(  0,\lambda\left(  \Gamma\right)
\right)  .$

The last step in our consideration is traditional \cite{Vas}: the even
extension of eigenfunctions in $\Pi^{+}\left(  \varepsilon\right)  $ through
the Dirichlet conditions onto the domain $\Pi\left(  \varepsilon\right)  $
becomes a trapped mode in the whole problem (\ref{1.5})-(\ref{1.7}).

\subsection{Preliminary description of results.\label{sect1.3}}

The operator formulation of problem (\ref{1.5})-(\ref{1.7}) given in Section
\ref{sect3.2} permits to deal with its spectrum within the spectral theory of
self-adjoint operators in Hilbert space. If $\lambda_{\varepsilon}%
\in\mathbb{C}$ is a complex number and $\lambda_{\varepsilon}\notin
\overline{\mathbb{R}}_{+}=\left[  0,+\infty\right)  ,$ then evidently, the
inhomogeneous problem (\ref{1.5})-(\ref{1.7}) with data in the Lebesgue spaces
$L^{2}\left(  \Pi\left(  \varepsilon\right)  \right)  $ and $L^{2}\left(
\partial\Pi\left(  \varepsilon\right)  \right)  $ admits a unique generalized
solution in the Sobolev space $H^{1}\left(  \Pi\left(  \varepsilon\right)
\right)  $ (see the integral identity (\ref{3.3}) and cf. \cite{Lad}). This
fact means that $\mathbb{C}\backslash\overline{\mathbb{R}}_{+}$ implies the
resolvent set of the operator $\mathcal{L}_{\varepsilon}$ of problem
(\ref{1.5})-(\ref{1.7}). In Section \ref{sect3.3} we show that the closed real
positive semi-axis is covered with the continuous spectrum of $\mathcal{L}%
_{\varepsilon}$ (Lemma \ref{Lemma3.1}).

Under the assumption
\begin{equation}
\Pi\left(  \varepsilon\right)  =\left\{  x:\left(  x_{1,}-x_{2,}x_{3}\right)
\in\Pi\left(  \varepsilon\right)  \right\}  , \label{1.8}%
\end{equation}
which requires the symmetry of domain (\ref{1.8}) with respect to the middle
plane $\left\{  x:x_{2}=0\right\}  $ of the canal (cf. Fig. \ref{f3}, a, where
the dotted line indicates the symmetry axis of the transverse cross-section of
the canal), we treat the restriction $\mathcal{L}_{\varepsilon}^{0}$ of the
operator $\mathcal{L}_{\varepsilon}$ onto the subspace
\begin{equation}
\mathcal{H}^{0}=\left\{  \Phi\in H^{1}\left(  \Pi\left(  \varepsilon\right)
\right)  :\Phi\text{ is odd in }x_{2}\right\}  \label{1.9}%
\end{equation}
and associate with the operator $\mathcal{L}_{\varepsilon}^{0}$ a problem
obtained from (\ref{1.5})-(\ref{1.7}) by restricting onto a half of the domain
$\Pi\left(  \varepsilon\right)  ,$ for definiteness on the right half
(\ref{1.10}), and supplied with the artificial boundary condition
\begin{equation}
\Phi_{\varepsilon}^{+}\left(  x\right)  =0,\text{ \ \ \ }x\in\varpi^{0}\left(
\varepsilon\right)  , \label{1.11}%
\end{equation}
on the artificially generated surface $\varpi^{0}\left(  \varepsilon\right)
=\left\{  x\in\Pi\left(  \varepsilon\right)  :x_{2}=0\right\}  .$ Such the
restricted problem is further referred as the problem (\ref{1.5})-(\ref{1.7}),
(\ref{1.11}) on the domain $\Pi^{+}\left(  \varepsilon\right)  .$

In Section \ref{sect3.3}, owing to the Dirichlet boundary conditions
(\ref{1.11}), we find out a threshold $\lambda\left(  \Gamma\right)  >0$,
depending only on the cross-section $\Gamma$, such that the continuous
spectrum of $\mathcal{L}_{\varepsilon}^{0}$ implies the ray $\left[
\lambda\left(  \Gamma\right)  ,+\infty\right)  \subset\mathbb{R}_{+}$ while
the segment $\left[  0,\lambda\left(  \Gamma\right)  \right)  $ contains only
the discrete spectrum of $\mathcal{L}_{\varepsilon}^{0}.$

Note that the odd extension $\Phi_{\varepsilon}$ of an eigenfunction
$\Phi_{\varepsilon}^{+}$ of the problem in $\Pi^{+}\left(  \varepsilon\right)
$ becomes an eigenfunction of the problem in $\Pi\left(  \varepsilon\right)  $
corresponding to the same eigenvalue $\lambda_{\varepsilon}=\lambda
_{\varepsilon}^{+}.$ Based on the above-mentioned observations, we prove in
Section \ref{sect3.4} the main result of the paper.

\begin{theorem}
\label{Theorem1.1}Under the geometrical assumptions (\ref{1.1}), (\ref{1.4})
and (\ref{1.8}), for any $d>0$ and $N\in\mathbb{N}:=\left\{  1,2,....\right\}
,$ there exists $\varepsilon\left(  d,N\right)  >0$ such that in the case
$\varepsilon\in\left(  0,\varepsilon\left(  d,N\right)  \right)  $ problem
(\ref{1.5})-(\ref{1.7}) has at least $N$ eigenvalues $\lambda_{\varepsilon
}^{\left(  1\right)  },...,\lambda_{\varepsilon}^{\left(  N\right)  }$ in the
interval $\left(  0,d\right)  \subset\mathbb{R}_{+}.$ The corresponding
eigenfunctions $\Phi_{\varepsilon}^{\left(  1\right)  },...,\Phi_{\varepsilon
}^{\left(  N\right)  }$ decay exponentially at infinity and, therefore, imply
so-called trapped modes in the linear theory of water-waves.
\end{theorem}

We emphasize that the eigenvalues in Theorem \ref{Theorem1.1} lie in the
continuous spectrum of the operator\ $\mathcal{L}_{\varepsilon}.$

Our approach does not require any other \textit{global} geometry assumption on
the shape of the body $\mathbf{\Theta}\left(  \varepsilon\right)  $ whilst the
symmetric cross-section $\Gamma$ of the canal is arbitrary. Moreover, any
given large number of trapped modes with the frequencies in any preadjusted
small interval can be obtained.

\section{Asymptotics of eigenvalues of the spectral problem in the thin
domain\label{sect2}}

\subsection{Formal asymptotic analysis.\label{sect2.1}}

\begin{figure}
[ptb]
\begin{center}
\includegraphics[
height=1.3578in,
width=1.3526in
]%
{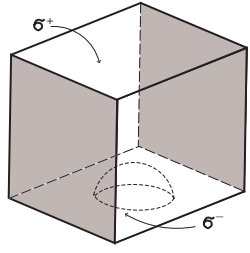}%
\caption{The periodicity cell $\Sigma.$}%
\label{f5}%
\end{center}
\end{figure}

We employ the standard asymptotic expansions of solutions in thin domains
(see, e.g., \cite{na192},\cite[Ch.7]{Nabook})%
\begin{equation}
\alpha_{\varepsilon}\sim\varepsilon\tau,\text{ \ }u_{\varepsilon}\left(
x\right)  \sim w\left(  y\right)  +\varepsilon w_{1}\left(  y,\xi\right)
+\varepsilon^{2}w_{2}\left(  y,\xi\right)  , \label{2.1}%
\end{equation}
where $\tau$ and $w,w_{j}$ are a number and functions to be determined and
$\xi$ stands for the "fast" variables%
\begin{equation}
\xi=\left(  \eta,\zeta\right)  ,\text{ \ }\eta=\varepsilon^{-1}y,\text{
\ }\zeta=\varepsilon^{-1}z. \label{2.2}%
\end{equation}
We insert formulae (\ref{2.1}) into equation (\ref{1.12}) and the boundary
conditions (\ref{1.13}) and gather coefficients on similar powers of the small
parameter $\varepsilon.$ Since the derivatives in $y_{i}$ and $z$ of the
function $\left(  y,z\right)  \mapsto W\left(  y,\varepsilon^{-1}%
y,\varepsilon^{-1}z\right)  $ are equal to
\[
\varepsilon^{-1}\frac{\partial W}{\partial\xi_{i}}\left(  y,\xi\right)
+\frac{\partial W}{\partial y_{i}}\left(  y,\xi\right)  \text{ \ and
\ }\varepsilon^{-1}\frac{\partial W}{\partial\zeta}\left(  y,\xi\right)  ,
\]
respectively, we obtain the following problems on the periodicity cell
$\Sigma$ with the parameter $y\in\omega:$%
\begin{equation}
-\Delta_{\xi}w\left(  y\right)  =0,\text{ \ }\xi\in\Sigma,\text{
\ \ \ }\partial_{n\left(  \xi\right)  }w\left(  y\right)  =0,\text{ \ }\xi
\in\sigma^{+}\cup\sigma^{-}; \label{2.3}%
\end{equation}%
\begin{align}
-\Delta_{\xi}w_{1}\left(  y,\xi\right)   &  =2\nabla_{\eta}\cdot\nabla
_{y}w\left(  y\right)  ,\text{ \ }\xi\in\Sigma,\label{2.4bis}\\
\partial_{n\left(  \xi\right)  }w_{1}\left(  y,\xi\right)   &  =-n^{\bullet
}\left(  \xi\right)  \cdot\nabla_{y}w\left(  y\right)  ,\text{ \ }\xi\in
\sigma^{+}\cup\sigma^{-};\nonumber
\end{align}%
\begin{align}
-\Delta_{\xi}w_{2}\left(  y,\xi\right)   &  =2\nabla_{\eta}\cdot\nabla
_{y}w_{1}\left(  y\right)  +\Delta_{y}w\left(  y\right)  ,\text{\ \ }\xi
\in\Sigma,\label{2.5}\\
\partial_{n\left(  \xi\right)  }w_{2}\left(  y,\xi\right)   &  =-n^{\bullet
}\left(  \xi\right)  \cdot\nabla_{y}w_{1}\left(  y\right)  ,\text{ \ }\xi
\in\sigma^{-},\nonumber\\
\partial_{\zeta}w_{2}\left(  y,\xi\right)   &  =\tau w\left(  y\right)
,\ \xi\in\sigma^{+}.\nonumber
\end{align}
Here $n=\left(  n_{1},n_{2,}n_{3}\right)  $ is the outward unit normal to the
upper $\sigma^{+}$ and lower $\sigma^{-\text{ }}$ bases of the cell $\Sigma$
(see Fig. \ref{f5} and compare with Fig. \ref{f4}), $n=\left(  0,0,1\right)  $
on $\sigma^{+}=\left\{  \xi:\eta\in\sigma,\zeta=0\right\}  ,$ and $n^{\bullet
}=\left(  n_{1},n_{2}\right)  $ so that $n^{\bullet}=\left(  0,0\right)  $ on
$\sigma^{+},$ $\nabla_{y}=\left(  \partial/\partial y_{1},\partial/\partial
y_{2}\right)  .$%

Problems (\ref{2.3})-(\ref{2.5}) are also supplied with the periodicity
conditions on the opposite lateral sides $\left\{  \xi\in\partial\Sigma
:\eta_{i}=\pm a_{i}/2\right\}  ,$ $i=1,2,$ of the cell (a couple of them is
overshadowed in Fig. \ref{f5}). Note that $\sigma^{-}$ is the surface which
completes these sides and the rectangular "cover" $\sigma^{+}$ up to the whole
boundary $\partial\Sigma.$ We do not write explicitly the periodicity
conditions but always deal with solutions which are $a_{i}-$periodic in the
variables $\eta_{i},$ $i=1,2.$

Equations (\ref{2.3}) hold true because $w$ does not depend on the fast
variables $\xi$ in (\ref{2.2}). Since, evidently,
\[
\int_{\sigma^{-}}n_{i}\left(  \xi\right)  ds_{\xi}=\int_{\partial\Sigma}%
n_{i}\left(  \xi\right)  ds_{\xi}=\int_{\Sigma}\frac{\partial1}{\partial
\xi_{i}}d\xi=0,\text{ \ }i=1,2,
\]
problem (\ref{2.4bis}) admits a solution in the form%
\begin{equation}
w_{1}\left(  y,\xi\right)  =-\underset{i=1}{\overset{2}{\sum}}W_{i}\left(
\xi\right)  \frac{\partial w}{\partial y_{i}}\left(  y\right)  , \label{2.6}%
\end{equation}
where $W_{1}$ and $W_{2}$ raise the standard asymptotic corrector in the
theory of homogenization (see, e.g., \cite{BLP, SP}). Namely, $W_{i}$ is a
(periodic in $\eta$) solution of the model problem
\begin{equation}
-\Delta_{\xi}W_{i}\left(  \xi\right)  =0,\text{ \ }\xi\in\Sigma,\text{
\ \ \ }\partial_{n\left(  \xi\right)  }W_{i}\left(  \xi\right)  =n_{i}\left(
\xi\right)  ,\text{ \ }\xi\in\sigma^{+}\cup\sigma^{-}. \label{2.7}%
\end{equation}
We emphasize that, by definition, $n_{1}=n_{2}=0$ on $\sigma^{+}$ and,
according to the assumed smoothness of the lower base of the cell, the
periodic functions $W_{i}$ are infinitely differentiable.

We now consider problem (\ref{2.5}). Note that the factor $\varepsilon$ in the
representation (\ref{2.1}) of the eigenvalue $\alpha_{\varepsilon}$ was
introduced to fulfil the goal: the main asymptotic term of the right-hand side
$\alpha_{\varepsilon}u_{\varepsilon}\left(  x\right)  $ in the spectral
boundary condition (\ref{1.13}) of the Steklov type comes into a problem for
the asymptotic term $\varepsilon^{2}w_{2}$ in the expansion for the
eigenfunction $u_{\varepsilon}.$

The compatibility condition in problem (\ref{2.5}) reads
\begin{equation}
0=2\int_{\Sigma}\nabla_{\eta}\cdot\nabla_{y}w_{1}\left(  y,\xi\right)
d\xi+\left\vert \Sigma\right\vert \Delta_{y}w\left(  y\right)  -\int
_{\sigma^{-}}n^{\bullet}\left(  \xi\right)  \cdot\nabla_{y}w_{1}\left(
y,\xi\right)  ds_{\xi}+\beta\left\vert \sigma\right\vert w\left(  y\right)  ,
\label{2.8}%
\end{equation}
where $\left\vert \Sigma\right\vert =$meas$_{3}\Sigma$\ is the volume of the
cell $\Sigma$ and $\left\vert \sigma\right\vert =a_{1}a_{2}$ the area of the
cover $\sigma^{+}.$ Owing to (\ref{2.6}), equality (\ref{2.8}) can be
rewritten in the form
\begin{equation}
B\left(  \nabla_{y}\right)  w\left(  y\right)  \ \text{:=}\ -\nabla_{y}\cdot
b\nabla_{y}w\left(  y\right)  =\tau\left\vert \sigma\right\vert w\left(
y\right)  ,\text{\ \ }y\in\omega. \label{2.10}%
\end{equation}
Here $b$ is a matrix of size $2\times2$ with the entries%
\begin{align}
b_{ik}  &  =-%
{\displaystyle\int_{\Sigma}}
\left(  \dfrac{\partial W_{k}}{\partial\eta_{i}}\left(  \xi\right)
+\dfrac{\partial Wi}{\partial\eta_{k}}\left(  \xi\right)  \right)
d\xi+\left\vert \Sigma\right\vert \delta_{i,k}+%
{\displaystyle\int_{\sigma^{-}}}
W_{k}\left(  \xi\right)  \partial_{n\left(  \xi\right)  }W_{i}\left(
\xi\right)  ds_{\xi}=\label{2.9}\\
&  =%
{\displaystyle\int_{\Sigma}}
\left(  \delta_{i,k}-\dfrac{\partial W_{k}}{\partial\eta_{i}}\left(
\xi\right)  -\dfrac{\partial W_{i}}{\partial\eta_{k}}\left(  \xi\right)
+\nabla_{\xi}W_{k}\left(  \xi\right)  \cdot\nabla_{\xi}W_{i}\left(
\xi\right)  \right)  d\xi=\nonumber\\
&  =\left(  \nabla_{\xi}\left(  \xi_{k}-W_{k}\right)  ,\nabla_{\xi}\left(
\xi_{i}-W_{i}\right)  \right)  _{\Sigma}.\nonumber
\end{align}
By $\left(  \cdot,\cdot\right)  _{\Sigma}$ is denoted the natural inner
product in the Lebesgue space $L^{2}\left(  \Sigma\right)  $. The vector
functions $\nabla_{\xi}\left(  \xi_{1}+W_{1}\right)  $ and $\nabla_{\xi
}\left(  \xi_{2}+W_{2}\right)  $ are linear independent because $W_{1}$ and
$W_{2}$ are periodic in $\eta.$ Thus, the matrix $b$ with entries (\ref{2.9})
implies a Gram matrix, i.e. it is positive definite and symmetric and,
therefore, $B\left(  \nabla_{y}\right)  $ is a second order elliptic
differential operator.

In order to satisfy the Dirichlet conditions (\ref{1.14}) on the lateral side
of the plate $\Omega_{\varepsilon}$, we subject the function $w$ in
(\ref{2.1}) to the boundary condition%
\begin{equation}
w\left(  y\right)  =0,\ y\in\partial\omega. \label{2.11}%
\end{equation}
We call (\ref{2.10}), (\ref{2.11}) the resultant spectral problem.

If $\beta$ and $w$ are an eigenvalue and the corresponding eigenfunction of
problem (\ref{2.10}), (\ref{2.11}), the compatibility condition (\ref{2.8}) is
met and problem (\ref{2.5}) admits a solution. This completes the asymptotic
expansion (\ref{2.1}).

\subsection{Spectrum of the resultant problem.}

Problem (\ref{2.10}), (\ref{2.11}) can be reformulated as the integral
identity \cite{Lad}%
\begin{equation}
\left(  b\nabla_{y}w,\nabla_{y}v\right)  _{\omega}=\tau\left\vert
\sigma\right\vert \left(  \omega,v\right)  _{\omega}\text{,} \label{2.12}%
\end{equation}
the left-hand side of which implies an inner product in the subspace
$\mathring{H}^{1}\left(  \omega,\partial\omega\right)  $ of functions $w\in
H^{1}\left(  \omega\right)  $ satisfying condition (\ref{2.11}). Owing to the
compact embedding $H^{1}\left(  \omega\right)  \subset L^{2}\left(
\omega\right)  ,$ spectrum of the operator, associated with the bi-linear form
$\left(  b\nabla_{y}\cdot,\nabla_{y}\cdot\right)  _{\omega}$ (see
\cite[\S 10.1]{BiSo}), is descrete and form the positive monotone unbounded
sequence
\begin{equation}
0<\tau^{\left(  1\right)  }<\tau^{\left(  2\right)  }\leq\tau^{\left(
3\right)  }\leq...\leq\tau^{\left(  p\right)  }\leq....\rightarrow
+\infty\label{2.13}%
\end{equation}
where eigenvalues are repeated according to their multiplicity.

The corresponding eigenfunctions $w^{\left(  1\right)  },$ $w^{\left(
2\right)  },...,$ $w^{\left(  p\right)  },...$ can be subject to the
normalization and orthogonality condition
\begin{equation}
\left(  b\nabla_{y}w^{\left(  p\right)  },\nabla_{y}w^{\left(  q\right)
}\right)  _{\omega}+\left\vert \sigma\right\vert \left(  w^{\left(  p\right)
},w^{\left(  q\right)  }\right)  _{\omega}=\delta_{p,q} \label{2.14}%
\end{equation}
where $p$, $q\in\mathbb{N}$ and $\delta_{p,q}$ is Kronecker's symbol. The
first eigenvalue $\tau^{\left(  1\right)  }$ is simple due to the strong
maximum principle.

An affine transform of the coordinate system $y=\left(  y_{1\text{,}}%
y_{2}\right)  $ turns the differential operator $B\left(  \nabla_{y}\right)  $
on the left of (\ref{2.10}) into the Laplace operator while the rectangle
$\omega$ becomes a parallelogram. A harmonic function, which has the finite
Dirichlet integral and vanishes at both sides of an angle with the opening
$\psi\in\left(  0,\pi\right)  ,$ possesses the worst singularity $Kr^{\pi
/\psi}\sin\left(  \pi\psi^{-1}\varphi\right)  $ where $\left(  r,\varphi
\right)  $ is the polar coordinate system and $K\in\mathbb{R}$ (see, e.g.,
\cite{Gris}, and introductory chapters in \cite{NaPL}, \cite{KoMaRo})$.$ Thus,
the theory of elliptic boundary value problems in domains with piecewise
smooth boundaries, especially, a result in \cite{MaPL2}, ensures the following assertion.

\begin{lemma}
\label{Lemma 2.1}The eigenfunction $w^{\left(  p\right)  }\in\mathring{H}%
^{1}\left(  \omega\right)  $ of problem (\ref{2.10}), (\ref{2.11}) verifies
the estimates%
\begin{equation}
\left\vert \nabla_{y}^{k}w^{\left(  p\right)  }\left(  y\right)  \right\vert
\leq c_{p,k}R\left(  x\right)  ^{1+\rho-k}\text{, \ }k\in\mathbb{N}%
_{0}=\left\{  0,1,2,...\right\}  , \label{2.99}%
\end{equation}
where $\rho\in\left(  0,1\right)  $ is a number depending on the matrix $b$
with entries (\ref{2.9}) and $R\left(  x\right)  $ is the distance from a
point $x\in\overline{\omega}$ to the nearest among the tops of the rectangle
$\omega$. In particular, $w^{\left(  p\right)  }$ belongs to the Sobolev space
$H^{2}\left(  \omega\right)  $ and the H\"{o}lder space $C^{1,\rho}\left(
\omega\right)  .$
\end{lemma}

Recall the definition of the H\"{o}lder norm
\begin{equation}
\left\Vert w;C^{l,\rho}\left(  \omega\right)  \right\Vert =\overset
{l}{\underset{k=0}{\sum}~}\underset{y\in\omega}{\sup}\left\vert \nabla_{y}%
^{k}w\left(  y\right)  \right\vert +\underset{y,\mathbf{y}\in\omega}{\sup
}\left\vert y-\mathbf{y}\right\vert ^{-\rho}\left\vert \nabla_{y}^{l}w\left(
y\right)  -\nabla_{\mathbf{y}}^{l}w\left(  \mathbf{y}\right)  \right\vert .
\label{2.97}%
\end{equation}

\subsection{Operator formulation of the problem in $\Omega_{\varepsilon.}$}

Aiming to justify asymptotic expansions constructed in Section \ref{sect2.1},
we endow the Sobolev space%
\begin{equation}
\mathring{H}^{1}\left(  \Omega_{\varepsilon},\mathbf{\Upsilon}_{\varepsilon
}\right)  =\left\{  u\in H^{1}\left(  \Omega_{\varepsilon}\right)  :u\left(
x\right)  =0,\text{ }x\in\mathbf{\Upsilon}_{\varepsilon}\right\}  \label{2.00}%
\end{equation}
with the specific inner product%
\begin{equation}
\left\langle u,v\right\rangle _{\varepsilon}=\left(  \nabla_{x}u,\nabla
_{x}v\right)  _{\Omega_{\varepsilon}}+\varepsilon\left(  u,v\right)  _{\omega
}. \label{2.17}%
\end{equation}
In the obtained Hilbert space $\mathcal{H}_{\Omega}^{\varepsilon}$ we
introduce the operator $\mathcal{B}_{\varepsilon}$ by the formula
\begin{equation}
\left\langle \mathcal{B}_{\varepsilon}u,v\right\rangle _{\varepsilon}=\left(
u,v\right)  _{\omega^{+}},\ \ u,v\in\mathcal{H}_{\Omega}^{\varepsilon}.
\label{2.16}%
\end{equation}
This operator is positive continuous and symmetric, therefore, self-adjoint.
It is compact due to the compactness of the embedding $H^{1}\left(
\Omega_{\varepsilon}\right)  \subset L^{2}\left(  \partial\Omega_{\varepsilon
}\right)  .$ The norm of $\mathcal{B}_{\varepsilon}$ is less than
$\varepsilon^{-1}.$ Thus, the spectrum of $\mathcal{B}_{\varepsilon}$ is
descrete and forms the positive infinitesimal sequence $\left\{
\beta_{\varepsilon}^{\left(  j\right)  }\right\}  _{j\in\mathbb{N}},$%
\begin{equation}
\varepsilon^{-1}>\beta_{\varepsilon}^{\left(  1\right)  }>\beta_{\varepsilon
}^{\left(  2\right)  }\geq\beta_{\varepsilon}^{\left(  3\right)  }\geq
...\beta_{\varepsilon}^{\left(  p\right)  }\geq....\rightarrow0^{+},
\label{2.17bis}%
\end{equation}
where eigenvalues are listed according to their multiplicity and again the
first eigenvalue is simple by virtue of the strong maximum principle. The
corresponding eigenfunctions $u_{\varepsilon}^{\left(  p\right)  }$ can be
subject to the normalization and orthogonality conditions%
\begin{equation}
\left\langle u_{\varepsilon}^{\left(  p\right)  },u_{\varepsilon}^{\left(
q\right)  }\right\rangle _{\varepsilon}=\delta_{p,q},\text{ \ }p,q\in
\mathbb{N}. \label{2.18}%
\end{equation}

\begin{remark}
\label{Remark2.2}The variational formulation of problem (\ref{1.12}%
)-(\ref{1.14})%
\begin{equation}
\left(  \nabla_{x}u_{\varepsilon},\nabla_{x}v\right)  _{\Omega_{\varepsilon}%
}=\alpha_{\varepsilon}\left(  u_{\varepsilon},v\right)  _{\omega^{_{+}}%
},\text{\ \ }v\in\mathring{H}^{1}\left(  \Omega_{\varepsilon};\mathbf{\Upsilon
}_{\varepsilon}\right)  , \label{2.19}%
\end{equation}
is equivalent to the abstract equation%
\begin{equation}
\mathcal{B}_{\varepsilon}u_{\varepsilon}=\beta_{\varepsilon}u_{\varepsilon}%
\in\mathcal{H}_{\Omega}^{\varepsilon} \label{2.20}%
\end{equation}
with the new spectral parameter%
\begin{equation}
\beta_{\varepsilon}=\left(  \alpha_{\varepsilon}+\varepsilon\right)  ^{-1}.
\label{2.21}%
\end{equation}
Formula (\ref{2.21}) relates only the descrete spectra (\ref{2.17}) and
(\ref{1.88}). Although the operator $\mathcal{B}_{\varepsilon}$ has the
infinite-dimensional kernel $\mathring{H}^{1}\left(  \Omega_{\varepsilon
};\mathbf{\Upsilon}_{\varepsilon}\cup\omega^{+}\right)  ,$ this kernel does
not influence the spectrum of problem (\ref{2.12}) because $\beta
=0\mapsto\alpha=\varepsilon-\beta^{-1}=\infty.$
\end{remark}

Justification of asymptotics is based in the next sections on the following
classical result known as the lemma on "almost eigenvalues and
eigenfunctions", a proof can be found in \cite{ViLu} and \cite{BiSo}.

\begin{lemma}
\label{Lemma 2.3} Let\textbf{\ }$\mathbf{u}\in\mathcal{H}_{\Omega
}^{\varepsilon}$ and $\mathbf{b}\in\mathbb{R}_{+}$ satisfy
\[
\left\Vert \mathbf{u};\mathcal{H}_{\Omega}^{\varepsilon}\right\Vert =1,\text{
\ }\left\Vert \mathcal{B}_{\varepsilon}\mathbf{u}-\mathbf{bu};\mathcal{H}%
_{\Omega}^{\varepsilon}\right\Vert =\delta<\mathbf{b}.
\]
Then at least one eigenvalue $\beta_{\varepsilon}^{\left(  q\right)  }$of the
operator $\mathcal{B}_{\varepsilon}$ verifies the inequality%
\[
\left\vert \beta_{\varepsilon}^{\left(  q\right)  }-\mathbf{b}\right\vert
\leq\delta.
\]
Moreover, for any $\delta_{1}\in\left(  \delta,\mathbf{b}\right)  $, there
exist coefficients $f_{p}$ such that
\[
\sum\left\vert f_{p}\right\vert ^{2}=1,\text{ }\left\Vert \mathbf{u}-\sum
f_{p}u_{\varepsilon}^{\left(  p\right)  };\mathcal{H}_{\Omega}^{\varepsilon
}\right\Vert \leq2\frac{\delta}{\delta_{1}}%
\]
where $\sum$ means summation over all eigenvalues of the operator
$\mathcal{B}_{\varepsilon}$ in the segment $\left[  \mathbf{b}-\delta
_{1},\mathbf{b}+\delta_{1}\right]  $ and $u_{\varepsilon}^{\left(  p\right)
}$ are corresponding eigenfunctions under condition (\ref{2.18})$.$
\end{lemma}

\subsection{Approximation solutions.\label{sect2.4}}

According to (\ref{2.1}) and (\ref{2.21}), we take
\begin{equation}
\mathbf{b}=\varepsilon^{-1}\left(  \tau^{\left(  k\right)  }+1\right)  ,\text{
\ }\mathbf{u}_{\varepsilon}^{\left(  p\right)  }\left(  x\right)  =\left\Vert
\mathbf{U}_{\varepsilon}^{\left(  p\right)  };\mathcal{H}_{\Omega
}^{\varepsilon}\right\Vert ^{-1}\mathbf{U}_{\varepsilon}^{\left(  p\right)
}\left(  x\right)  , \label{2.22}%
\end{equation}%
\begin{equation}
\mathbf{U}_{\varepsilon}^{\left(  p\right)  }\left(  x\right)  =w^{\left(
p\right)  }\left(  y\right)  +\varepsilon X_{\varepsilon}\left(  y\right)
U_{\varepsilon}^{\left(  p\right)  }\left(  x\right)  ,\text{ \ }%
U_{\varepsilon}^{\left(  p\right)  }\left(  x\right)  =\overset{2}%
{\underset{i=1}{\sum}}W_{i}\left(  \varepsilon^{-1}x\right)  \frac{\partial
w^{\left(  p\right)  }}{\partial y_{i}}\left(  y\right)  \label{2.23}%
\end{equation}
as an approximate solution of the spectral abstract equation (\ref{2.20}). In
(\ref{2.22}) $\tau^{\left(  k\right)  }$ is an eigenvalue of the resultant
problem with multiplicity $\varkappa_{k},$ i.e.,%
\begin{equation}
\tau^{\left(  k-1\right)  }<\tau^{\left(  k\right)  }=...=\tau^{\left(
k+\varkappa_{k}-1\right)  }<\tau^{\left(  k+\varkappa_{k}\right)  }
\label{2.24}%
\end{equation}
in the sequence (\ref{2.13}), $X_{\varepsilon}$ is a smooth cut-off function
on $\omega$ which is equal to $1$ outside the $\varepsilon-$neighborhood of
$\partial\omega$ and vanishes in the vicinity of $\partial\omega,$ $\left\vert
\nabla_{y}^{j}X_{\varepsilon}\left(  y\right)  \right\vert \leq c\varepsilon
^{-j},$ e.g.,%
\begin{equation}
X_{\varepsilon}\left(  x\right)  =X_{\varepsilon}^{1}\left(  x_{1}\right)
X_{\varepsilon}^{2}\left(  x_{2}\right)  ,\ \ X_{\varepsilon}^{i}\left(
x_{i}\right)  =\left\{
\begin{tabular}
[c]{l}%
$1$ \ \ \ for $\left\vert x_{i}\right\vert <\tfrac{1}{2}a_{i}-\varepsilon,$\\
$0$ $\ \ \ $for $\left\vert x_{i}\right\vert >\tfrac{1}{2}a_{i}-\tfrac{1}%
{2}\varepsilon.$%
\end{tabular}
\ \ \right.  \label{2.25}%
\end{equation}
Furthermore, $p=k,....k+\varkappa_{k}-1$ and $w^{\left(  k\right)
},...,w^{\left(  k+\varkappa_{k}-1\right)  }$ are eigenfunctions of problem
(\ref{2.10}), (\ref{2.11}) corresponding to $\tau^{\left(  k\right)  }$ and
verifying conditions (\ref{2.14}). In other words, formulae (\ref{2.22}),
(\ref{2.23}) deliver $\varkappa_{k}$ different approximation solutions of
(\ref{2.20}).

We proceed with calculation of the inner products $\left\langle \mathbf{U}%
_{\varepsilon}^{\left(  p\right)  },\mathbf{U}_{\varepsilon}^{\left(
q\right)  }\right\rangle _{\varepsilon};$ here and in the sequel
$p,q=k,...,k+\varkappa_{k}-1.$ Since $w^{\left(  p\right)  }\in H^{2}\left(
\omega\right)  ,$ $W^{i}\in C^{1}\left(  \Sigma\right)  $ and
\[
\nabla_{y}U_{\varepsilon}^{\left(  p\right)  }\left(  x\right)  =\overset
{2}{\underset{i=1}{\sum}}\left(  \varepsilon^{-1}\nabla_{\eta}W_{i}\left(
\xi\right)  \frac{\partial w^{\left(  p\right)  }}{\partial y_{i}}\left(
y\right)  +W_{i}\left(  \xi\right)  \nabla_{y}\frac{\partial w^{\left(
p\right)  }}{\partial y_{i}}\left(  y\right)  \right)  ,
\]
we readily obtain%
\[%
\begin{tabular}
[c]{l}%
$\left\Vert w^{\left(  p\right)  };\mathcal{H}_{\Omega}^{\varepsilon
}\right\Vert \leq c\varepsilon^{1/2},$\\
\\
$\left\Vert U^{\left(  p\right)  };L^{2}\left(  \Omega_{\varepsilon}\right)
\right\Vert +\varepsilon^{1/2}\left\Vert U^{\left(  p\right)  };L^{2}\left(
\omega_{+}\right)  \right\Vert +\varepsilon^{1/2}\left\Vert \nabla
_{x}U^{\left(  p\right)  };L^{2}\left(  \Omega_{\varepsilon}\right)
\right\Vert \leq c\varepsilon^{1/2}.$%
\end{tabular}
\ \ \ \ \ \
\]
Moreover,%
\[
\left\Vert U^{\left(  p\right)  }\nabla_{y}X_{\varepsilon};L^{2}\left(
\Omega_{\varepsilon}\right)  \right\Vert ^{2}\leq c\varepsilon^{-2}%
\text{meas}_{3}\left\{  x\in\Omega_{\varepsilon}:\text{dist}\left(
y,\partial\omega\right)  \leq c\varepsilon\right\}  \leq c
\]
and analogously%
\[
\left\Vert \left(  1-X_{\varepsilon}\right)  U^{\left(  p\right)  };\text{
}L^{2}\left(  \Omega_{\varepsilon}\right)  \right\Vert ^{2}\leq c\varepsilon
^{2}.
\]
These inequalities allow to estimate directly certain terms on the right-hand
side of the equality
\[%
\begin{tabular}
[c]{c}%
$\left\langle \mathbf{U}_{\varepsilon}^{\left(  p\right)  },\mathbf{U}%
_{\varepsilon}^{\left(  q\right)  }\right\rangle _{\varepsilon}=\left(
\nabla_{y}w^{\left(  p\right)  }+\varepsilon U^{\left(  p\right)  }\nabla
_{y}X_{\varepsilon}+\varepsilon X_{\varepsilon}\nabla_{y}U^{\left(  p\right)
},\text{ }\nabla_{y}w^{\left(  q\right)  }+\left.  +\varepsilon U^{\left(
q\right)  }\nabla_{y}X_{\varepsilon}+\varepsilon X_{\varepsilon}\nabla
_{y}U^{\left(  q\right)  }\right)  _{\Omega_{\varepsilon}}+\right.  $\\
$+\left(  \varepsilon X_{\varepsilon}\partial_{z}U^{\left(  p\right)
},\varepsilon X_{\varepsilon}\partial_{z}U^{\left(  q\right)  }\right)
_{\Omega_{\varepsilon}}+\varepsilon\left(  w^{\left(  p\right)  }+\varepsilon
X_{\varepsilon}U^{\left(  p\right)  },w^{\left(  q\right)  }+\varepsilon
X_{\varepsilon}U^{\left(  q\right)  }\right)  _{\omega_{+}}$%
\end{tabular}
\ \ \ \ \ \ \ \
\]
and conclude that
\begin{equation}
\left\vert \left\langle \mathbf{U}_{\varepsilon}^{\left(  p\right)
},\mathbf{U}_{\varepsilon}^{\left(  q\right)  }\right\rangle _{\varepsilon
}-J_{pq}-\varepsilon\left(  w^{\left(  p\right)  },w^{\left(  q\right)
}\right)  _{\omega}\right\vert \leq c\varepsilon^{3/2}. \label{2.26}%
\end{equation}
The formula%
\begin{equation}
\left\vert J_{pq}-\varepsilon\left\vert \sigma\right\vert ^{-1}\left(
B\nabla_{y}w^{\left(  p\right)  },\nabla_{y}w^{\left(  q\right)  }\right)
_{\omega}\right\vert \leq c\varepsilon^{1+\min\left\{  \rho,1/2\right\}  }
\label{2.27}%
\end{equation}
for the integral
\[
J_{pq}=\left(  \nabla_{y}w^{\left(  p\right)  }+\underset{i=1}{\overset
{2}{\sum}}\frac{\partial w^{\left(  p\right)  }}{\partial y_{i}}\nabla_{\xi
}W_{i},\nabla_{y}w^{\left(  q\right)  }+\underset{i=1}{\overset{2}{\sum}}%
\frac{\partial w^{\left(  q\right)  }}{\partial y_{k}}\nabla_{\xi}%
W_{k}\right)  _{\Omega_{\varepsilon}}%
\]
follows from the next lemma where it is necessary to put
\[
Z\left(  \xi\right)  =\nabla_{\xi}\left(  \xi_{i}+W_{i}\left(  \xi\right)
\right)  \cdot\nabla_{\xi}\left(  \xi_{k}+W_{k}\left(  \xi\right)  \right)
,\text{ \ }Y\left(  y\right)  =\frac{\partial w^{\left(  p\right)  }}{\partial
y_{i}}\left(  y\right)  \frac{\partial w^{\left(  q\right)  }}{\partial y_{k}%
}\left(  y\right)  .
\]
Note that $\rho$ is the exponent in Lemma \ref{Lemma 2.1} and formula
(\ref{2.9}) is used to detect the subtrahend on the left of (\ref{2.27}). The
following result is known (cf. \cite{BLP, SP}) so that we only adapt a
standard proof for the H\"{o}lder continuous multiplier $Y$ in the integrand.

\begin{lemma}
\label{Lemma2.4} Let $Z\in L^{\infty}\left(  \Sigma\right)  $ and $Y\in
C^{0,\rho}\left(  \omega\right)  ,$%
\begin{equation}
\overline{Z}=\left\vert \Sigma\right\vert ^{-1}%
{\displaystyle\int_{\Sigma}}
Z\left(  \xi\right)  d\xi. \label{2.28}%
\end{equation}
Then%
\begin{equation}
\left\vert \int_{\Omega_{\varepsilon}}Z\left(  \frac{x}{\varepsilon}\right)
Y\left(  y\right)  dx-\varepsilon\frac{\left\vert \Sigma\right\vert
}{\left\vert \sigma\right\vert }\overline{Z}\int_{\omega}Y\left(  y\right)
dy\right\vert \leq c\varepsilon^{1+\rho}. \label{2.29}%
\end{equation}

\end{lemma}

\textbf{Proof.} According to (\ref{1.4}), we have%
\begin{align*}%
{\displaystyle\int_{\Omega_{\varepsilon}}}
Z\left(  \dfrac{x}{\varepsilon}\right)  Y\left(  y\right)  dx  &
=\sum\limits_{\nu:\left\vert \nu_{i}\right\vert \leq N_{i}}%
{\displaystyle\int_{\Sigma_{\varepsilon}^{\nu}}}
Z\left(  \dfrac{x}{\varepsilon}\right)  Y\left(  y\right)  dx=\\
&  =\underset{\nu:\left\vert \nu_{i}\right\vert \leq N_{i}}{\sum}%
{\displaystyle\int_{\Sigma_{\varepsilon}^{\nu}}}
Z\left(  \dfrac{x}{\varepsilon}\right)  dx\left(  Y\left(  y^{\nu}\right)
+O\left(  \varepsilon^{\rho}\right)  \right)  =\varepsilon^{3}\left\vert
\Sigma\right\vert \overline{Z}\sum\limits_{\nu:\left\vert \nu_{i}\right\vert
\leq N_{i}}\left(  Y\left(  y^{\nu}\right)  +O\left(  \varepsilon^{\rho
}\right)  \right)  =\\
&  =\varepsilon\dfrac{\left\vert \Sigma\right\vert }{\left\vert \sigma
\right\vert }\overline{Z}\sum\limits_{\nu:\left\vert \nu_{i}\right\vert \leq
N_{i}}\left\vert \sigma_{\varepsilon}^{+\nu}\right\vert Y\left(  y^{\nu
}\right)  +O\left(  \varepsilon^{1+\rho}\right)  =\\
&  =\varepsilon\dfrac{\left\vert \Sigma\right\vert }{\left\vert \sigma
\right\vert }\overline{Z}\sum\limits_{\nu:\left\vert \nu_{i}\right\vert \leq
N_{i}}\left(
{\displaystyle\int_{\overset{+}{\sigma}_{\varepsilon}^{\nu}}}
Y\left(  y\right)  dy+O\left(  \varepsilon^{2+\rho}\right)  +O\left(
\varepsilon^{1+\rho}\right)  \right)  =\varepsilon\dfrac{\left\vert
\Sigma\right\vert }{\left\vert \sigma\right\vert }\overline{Z}%
{\displaystyle\int_{\omega}}
Y\left(  y\right)  dy+O\left(  \varepsilon^{1+\rho}\right)  .
\end{align*}
Here $\sigma_{\varepsilon}^{+\nu}=\left\{  y:\left\vert \varepsilon^{-1}%
y_{i}-\nu_{i}a_{i}\right\vert \leq a_{i}/2,\text{ }i=1,2\right\}  $ is the
cover of the cell (\ref{1.3}) and $y^{\nu}$ its mass center. The number of the
cells $\Sigma_{\varepsilon}^{\nu}$ composing the plate $\Omega_{\varepsilon}$
is less than $C\varepsilon^{-2}.$ Furthermore, we have used twice the
relation
\[
\left\vert Y\left(  y^{\nu}\right)  -Y\left(  y\right)  \right\vert \leq
c\varepsilon^{\rho}\text{ , \ \ }y\in\sigma_{i}^{\nu},
\]
inherited from the inclusion $Y\in C^{0,\rho}\left(  \omega\right)  $ and the
definition of the H\"{o}lder norm (\ref{2.97})$.\blacksquare$

Now formulae (\ref{2.26}), (\ref{2.29}) and (\ref{2.14}) ensure that
\begin{equation}
\left\vert \left\langle \mathbf{U}_{\varepsilon}^{\left(  p\right)
},\mathbf{U}_{\varepsilon}^{\left(  q\right)  }\right\rangle _{\varepsilon
}-\varepsilon\left\vert \sigma\right\vert ^{-1}\delta_{p,q}\right\vert \leq
c\varepsilon^{1+\min\left\{  \rho,1/2\right\}  }. \label{2.30}%
\end{equation}

\subsection{Calculating the discrepancy $\delta.$\label{sect2.5}}

According to (\ref{2.30}), we obtain
\begin{equation}
\left\Vert \mathbf{U}_{\varepsilon}^{\left(  p\right)  };\mathcal{H}_{\Omega
}^{\varepsilon}\right\Vert \geq\frac{1}{2}\varepsilon^{1/2}\left\vert
\sigma\right\vert ^{-1/2} \label{2.31}%
\end{equation}
for a small $\varepsilon>0.$ Thus, by virtue of (\ref{2.22}), (\ref{2.23}),
(\ref{2.16}), (\ref{2.17}), we have
\begin{align}
\delta &  =\left\Vert \mathcal{B}_{\varepsilon}\mathbf{u}_{\varepsilon
}^{\left(  p\right)  }-\mathbf{bu}_{\varepsilon}^{\left(  p\right)
};\mathcal{H}_{\Omega}^{\varepsilon}\right\Vert =\left\Vert U_{\varepsilon
}^{\left(  p\right)  };\mathcal{H}_{\Omega}^{\varepsilon}\right\Vert
^{-1}\mathbf{b}\left\Vert \mathbf{b}^{-1}\mathcal{B}_{\varepsilon
}U_{\varepsilon}^{\left(  p\right)  }-U_{\varepsilon}^{\left(  p\right)
};\mathcal{H}_{\Omega}^{\varepsilon}\right\Vert =\label{2.32}\\
&  =\left\Vert U_{\varepsilon}^{\left(  p\right)  };\mathcal{H}_{\Omega
}^{\varepsilon}\right\Vert ^{-1}\mathbf{b~}\sup\left\vert \varepsilon\left(
\tau^{\left(  k\right)  }+1\right)  \left\langle \mathcal{B}_{\varepsilon
}U_{\varepsilon}^{\left(  p\right)  },V\right\rangle _{\varepsilon
}-\left\langle U_{\varepsilon}^{\left(  p\right)  },V\right\rangle
_{\varepsilon}\right\vert \leq\nonumber\\
&  \leq2\varepsilon^{-3/2}\left\vert \sigma\right\vert ^{-1/2}\left(
\tau^{\left(  k\right)  }+1\right)  ^{-1}\sup\left\vert \left\langle
\mathcal{\nabla}_{x}U_{\varepsilon}^{\left(  p\right)  },\nabla_{x}%
V\right\rangle _{\Omega_{\varepsilon}}+\varepsilon\tau^{\left(  k\right)
}\left(  U_{\varepsilon}^{\left(  p\right)  },V\right)  _{\omega}\right\vert
,\nonumber
\end{align}
where the supremum is calculated over all $V\in\mathcal{H}_{\Omega
}^{\varepsilon}$ such that $\left\Vert V;\mathcal{H}_{\Omega}^{\varepsilon
}\right\Vert =1$. Furthermore,%
\begin{align}
I^{\left(  p\right)  }  &  =\left\langle \mathcal{\nabla}_{x}U_{\varepsilon
}^{\left(  p\right)  },\nabla_{x}V\right\rangle _{\Omega_{\varepsilon}%
}+\varepsilon\tau^{\left(  k\right)  }\left(  U_{\varepsilon}^{\left(
p\right)  },V\right)  _{\omega_{+}}=\label{2.33}\\
&  =-\left(  \Delta_{x}U_{\varepsilon}^{\left(  p\right)  },V\right)
_{\Omega_{\varepsilon}}+\left(  \mathcal{\partial}_{z}U_{\varepsilon}^{\left(
p\right)  },V\right)  _{\omega_{+}}+\varepsilon\tau^{\left(  k\right)
}\left(  U_{\varepsilon}^{\left(  p\right)  },V\right)  _{\omega_{+}}+\left(
\mathcal{\partial}_{n}U_{\varepsilon}^{\left(  p\right)  },V\right)
_{\omega_{-}\left(  \varepsilon\right)  }.\nonumber
\end{align}
To examine this expression we need auxiliary inequalities.

\begin{lemma}
\begin{enumerate}
\item \label{Lemma 2.5.1} Let $V\in\mathcal{H}_{\Omega}^{\varepsilon}$ and
\begin{equation}
\overline{V}\left(  \mathbf{y}\right)  =\left\vert \Sigma_{\varepsilon
}\right\vert ^{-1}\int_{\Sigma_{\varepsilon}\left(  \mathbf{y}\right)
}V\left(  x\right)  dx \label{2.34}%
\end{equation}
\noindent where
\[
\Sigma_{\varepsilon}\left(  \mathbf{y}\right)  =\left\{  x=\left(  y,z\right)
\in\Sigma_{\varepsilon}^{\infty}:\left\vert \mathbf{y}_{i}-y_{i}\right\vert
\leq\varepsilon a_{i}/2,\text{ \ \ }i=1,2\right\}
\]
\noindent and $V$ is extended by zero from $\Omega_{\varepsilon}$ onto the
thin periodic infinite layer. Then the inequality%
\begin{equation}
\left\Vert R_{\varepsilon}^{-1}\overline{V};L^{2}\left(  \omega\right)
\right\Vert +\left\Vert \nabla_{y}\overline{V};L^{2}\left(  \omega\right)
\right\Vert \leq c\varepsilon^{-1/2}\left\Vert V;\mathcal{H}_{\Omega
}^{\varepsilon}\right\Vert \label{2.35}%
\end{equation}
\noindent holds where $R_{\varepsilon}\left(  y\right)  =\varepsilon
+$dist$\left(  y,\partial\omega\right)  .$ Moreover,
\begin{equation}
\varepsilon^{-1}\left\Vert V-\overline{V};L^{2}\left(  \Omega_{\varepsilon
}\right)  \right\Vert +\varepsilon^{-1/2}\left\Vert V-\overline{V}%
;L^{2}\left(  \omega_{+}\cup\omega_{-}\left(  \varepsilon\right)  \right)
\right\Vert \leq c\left\Vert V;\mathcal{H}_{\Omega}^{\varepsilon}\right\Vert .
\label{2.36}%
\end{equation}

\item A function $V\in\mathcal{H}_{\Omega}^{\varepsilon}$ meets the relation%
\begin{equation}
\left\Vert R_{\varepsilon}^{-1}V;L^{2}\left(  \Omega_{\varepsilon}\right)
\right\Vert +\varepsilon^{1/2}\left\Vert R_{\varepsilon}^{-1}V;L^{2}\left(
\omega_{+}\cup\omega_{-}\left(  \varepsilon\right)  \right)  \right\Vert \leq
c\left\Vert V;\mathcal{H}_{\Omega}^{\varepsilon}\right\Vert . \label{2.37}%
\end{equation}
\noindent Here all constants depend on neither $V,$ nor $\varepsilon\in\left(
0,1\right]  .$
\end{enumerate}
\end{lemma}

\textbf{Proof.} First of all, we have%
\begin{align}%
{\displaystyle\int_{\omega}}
\left\vert \overline{V}\left(  \mathbf{y}\right)  \right\vert ^{2}d\mathbf{y}
&  \mathbf{\leq c}\varepsilon^{-6}%
{\displaystyle\int_{\omega}}
\left\vert
{\displaystyle\int_{\Sigma_{\varepsilon}\left(  \mathbf{y}\right)  }}
\left\vert V\left(  x\right)  \right\vert ^{2}dx\right\vert ^{2}%
d\mathbf{y\leq}\label{2.38}\\
&  \leq c\varepsilon^{-3}%
{\displaystyle\int_{\omega}}
{\displaystyle\int_{\Sigma_{\varepsilon}\left(  \mathbf{y}\right)  }}
\left\vert V\left(  y,z\right)  \right\vert ^{2}dydzd\mathbf{y}\leq\nonumber\\
&  \leq c\varepsilon^{-3}%
{\displaystyle\int_{\Omega_{\varepsilon}}}
{\displaystyle\int_{\sigma_{\varepsilon}\left(  y\right)  }}
d\mathbf{y}\left\vert V\left(  y,z\right)  \right\vert ^{2}dydz\leq
c\varepsilon^{-1}\left\Vert V;L^{2}\left(  \Omega_{\varepsilon}\right)
\right\Vert ^{2},\nonumber
\end{align}
where
\[
\sigma_{\varepsilon}\left(  y\right)  =\left\{  \mathbf{y:}\left\vert
\mathbf{y}-y_{i}\right\vert <\varepsilon a_{i}/2,\text{ \ }i=1,2\right\}  .
\]
Second,%
\begin{align*}
\left\vert \dfrac{\partial\overline{V}}{\partial\mathbf{y}_{i}}\left(
\mathbf{y}\right)  \right\vert  &  =\left\vert \Sigma_{\varepsilon}\right\vert
^{-1}\left\vert
{\displaystyle\int_{\sigma_{\varepsilon}^{i+}\left(  \mathbf{y}\right)  }}
V\left(  x\right)  ds_{x}-%
{\displaystyle\int_{\sigma_{\varepsilon}^{i-}\left(  \mathbf{y}\right)  }}
V\left(  x\right)  ds_{x}\right\vert \leq\\
&  \leq c\varepsilon^{-3}%
{\displaystyle\int_{\Sigma_{\varepsilon}\left(  \mathbf{y}\right)  }}
\left\vert \nabla_{x}V\left(  x\right)  \right\vert
dx\ \ \ \ \ \ \ \text{for\ almost\ all }\mathbf{y,}%
\end{align*}
where $\sigma_{\varepsilon}^{i\pm}\left(  \mathbf{y}\right)  =\left\{
x\in\partial\Sigma_{\varepsilon}\left(  \mathbf{y}\right)  :y_{i}%
=\mathbf{y}_{i}\pm\varepsilon a_{i}/2\right\}  $ are the opposite lateral
faces of the periodicity cell $\Sigma_{\varepsilon}\left(  \mathbf{y}\right)
$ (see Fig. \ref{f5})$.$ Now repeating calculation (\ref{2.38}) yields the
estimate of $\left\Vert \nabla_{y}\overline{V};L^{2}\left(  \Omega
_{\varepsilon}\right)  \right\Vert $ in (\ref{2.35}).

The support of function (\ref{2.34}) lies in the rectangle $\left\{
y:\left\vert y_{i}\right\vert \leq a_{i}\left(  2\varepsilon+1\right)
/2,\text{ \ }i=1,2\right\}  .$ Thus, integrating the one-dimensional Hardy
inequality%
\[%
{\displaystyle\int_{0}^{\infty}}
t^{-2}\left\vert \mathbf{V}\left(  t\right)  \right\vert ^{2}dt\leq4%
{\displaystyle\int_{0}^{\infty}}
\left\vert \frac{d\mathbf{V}}{dt}\left(  t\right)  \right\vert ^{2}dt,\text{
}\mathbf{V\in}C^{1}\left[  0,\infty\right)  ,\text{ }\mathbf{V}\left(
0\right)  =0,
\]
and using the completion argument bring the necessary estimate of the first
norm in (\ref{2.35})$.$

Dealing with the first norm in (\ref{2.36}), we compute%
\begin{equation}%
\begin{tabular}
[c]{l}%
$%
{\displaystyle\int_{\Omega_{\varepsilon}}}
\left\vert V\left(  y,z\right)  -\dfrac{1}{\left\vert \Sigma_{\varepsilon
}\right\vert }%
{\displaystyle\int_{\Sigma_{\varepsilon}\left(  y\right)  }}
V\left(  \mathbf{y},\mathbf{z}\right)  d\mathbf{y}d\mathbf{z}\right\vert
^{2}dydz=$\\
$\ \ \ \ \ \ \ =\dfrac{1}{\left\vert \Sigma_{\varepsilon}\right\vert ^{2}}%
{\displaystyle\int_{\Omega_{\varepsilon}}}
\left\vert
{\displaystyle\int_{\Sigma_{\varepsilon}\left(  y\right)  }}
\left(  V\left(  y,z\right)  -V\left(  \mathbf{y},\mathbf{z}\right)  \right)
d\mathbf{y}d\mathbf{z}\right\vert ^{2}dydz\leq$\\
$\ \ \ \ \ \ \leq c\varepsilon^{-3}%
{\displaystyle\int_{\Omega_{\varepsilon}}}
\varepsilon^{2}%
{\displaystyle\int_{\Sigma_{\varepsilon}\left(  y\right)  }}
\left\vert \nabla_{\mathbf{x}}V\left(  \mathbf{x}\right)  \right\vert
^{2}d\mathbf{x}dydz\leq c\varepsilon^{2}%
{\displaystyle\int_{\Omega_{\varepsilon}}}
\left\vert \nabla_{\mathbf{x}}V\left(  \mathbf{x}\right)  \right\vert
^{2}d\mathbf{x.}$%
\end{tabular}
\ \ \ \ \label{2.40}%
\end{equation}
For the second norm in (\ref{2.36}), we need to replace in (\ref{2.40}) the
integration set $\Omega_{\varepsilon}$ by $\omega_{+}$ and $\omega_{-}\left(
\varepsilon\right)  .$ As a result, the bound changes for $c\varepsilon
\left\Vert \nabla_{x}V;L^{2}\left(  \Omega_{\varepsilon}\right)  \right\Vert
^{2}.$

Inequality (\ref{2.37}) is a direct consequence of estimates (\ref{2.35}),
(\ref{2.36}) together with the evident relations $R_{\varepsilon}\left(
x\right)  ^{-1}\leq\varepsilon^{-1},$ $\varepsilon^{1/2}R_{\varepsilon}\left(
x\right)  ^{-1}\leq\varepsilon^{-1/2}$ and
\[
\varepsilon^{-1/2}\left\Vert R_{\varepsilon}^{-1}\overline{V};L^{2}\left(
\Omega_{\varepsilon}\right)  \right\Vert +\left\Vert R_{\varepsilon}%
^{-1}\overline{V};L^{2}\left(  \omega_{+}\cup\omega_{-}\left(  \varepsilon
\right)  \right)  \right\Vert \leq c\left\Vert R_{\varepsilon}^{-1}%
\overline{V};L^{2}\left(  \omega\right)  \right\Vert .\text{
\ \ \ \ \ \ \ \ \ \ \ \ }\blacksquare
\]
$\ \ \ \ \ \ \ \ \ \ \ \ \ \ \ \ \ \ \ \ \ \ \ \ \ \ \ \ \ \ \ \ \ \ \ \ \ \ \ \ \ \ \ \ \ \ \ \ \ \ \ \ \ \ \ \ \ \ \ \ \ \ \ \ \ \ \ \ \ \ \ \ \ \ \ \ \ \ \ \ \ \ \ \ \ \ \ \ \ \ \ \ \ \ \ \ \ \ \ \ \ \ \ \ \ \ \ \ \ \ \ \ \ \ \ \ \ \ \ \ \ \ \ \ \ \ \ \ \ \ \ \ \ \ \ \ \ \ \ \ \ \ \ \ \ \ \ \ \ \ \ \ \ \ \ \ \ \ \ \ \ \ \ \ \ \ \ \ \ \ \ \ \ \ \ \ \ \ \ \ \ \ \ \ \ \ \ \ \ \ \ \ \ \ \ \ \ \ \ \ \ \ \ \ \ \ \ \ \ \ \ \ \ \ \ \ \ \ \ \ \ \ \ \ \ \ \ \ \ \ \ \ \ \ \ \ \ \ \ \ \ \ $%

Now we are in position to simplify expression (\ref{2.33}) and, neglecting
inessential terms and changing $V$ for $\overline{V}$, to derive that%
\begin{equation}%
\begin{tabular}
[c]{l}%
$\left\vert I^{\left(  p\right)  }-\left(  \Delta_{y}w^{\left(  p\right)
}+2\overset{2}{\underset{i=1}{\sum}}\nabla_{\xi}W_{i}\cdot\nabla_{y}%
\dfrac{\partial w^{\left(  p\right)  }}{\partial y_{i}},X_{\varepsilon
}\overline{V}\right)  _{\Omega_{\varepsilon}}-\right.  $\\
$\ \ \ \ \ \ \left.  -\varepsilon\underset{i=1}{\overset{2}{\sum}}\left(
W_{i}n^{\bullet}\cdot\nabla_{y}\dfrac{\partial w^{\left(  p\right)  }%
}{\partial y_{i}},X_{\varepsilon}\overline{V}\right)  _{\omega^{-}\left(
\varepsilon\right)  }-\varepsilon\tau^{\left(  k\right)  }\left(  w^{\left(
p\right)  },X_{\varepsilon}\overline{V}\right)  _{\varpi}\right\vert \leq
c\varepsilon^{\rho+1/2}.$%
\end{tabular}
\ \ \ \ \ \label{2.39}%
\end{equation}
We start with the simplest term
\[
\left(  \mathcal{\partial}_{z}U_{\varepsilon}^{\left(  p\right)  }%
+\varepsilon\tau^{\left(  k\right)  }U_{\varepsilon}^{\left(  p\right)
},V\right)  _{\omega_{+}}=\left(  X_{\varepsilon}\underset{i=1}{\overset
{2}{\sum}}\dfrac{\partial w^{\left(  p\right)  }}{\partial y_{i}}\left(
\partial_{\zeta}W_{i}+\varepsilon^{2}\tau^{\left(  k\right)  }W_{i}\right)
+\varepsilon\tau^{\left(  k\right)  }w^{\left(  p\right)  },V\right)
_{\omega_{+}}=:I_{2}^{\left(  p\right)  }.
\]
Owing to (\ref{2.7}), we here have $\partial_{\zeta}W_{i}=0$ at $\zeta=0,$
and, hence,%
\begin{equation}%
\begin{tabular}
[c]{l}%
$\left\vert I_{2}^{\left(  p\right)  }-\varepsilon\tau^{\left(  k\right)
}\left(  w^{\left(  p\right)  },X_{\varepsilon}\overline{V}\right)  _{\varpi
}\right\vert \leq$\\
$\ \ \ \ \leq c\varepsilon\tau^{\left(  k\right)  }\left\vert \varepsilon
\left(  X_{\varepsilon}\underset{i=1}{\overset{2}{\sum}}\dfrac{\partial
w^{\left(  p\right)  }}{\partial y_{i}}W_{i},V\right)  _{\omega_{+}}+\left(
\left(  1-X_{\varepsilon}\right)  w^{\left(  p\right)  },V\right)
_{\omega_{+}}+\left(  X_{\varepsilon}w^{\left(  p\right)  },V-\overline
{V}\right)  _{\omega_{+}}\right\vert \leq$\\
$\ \ \ \ \leq c\varepsilon\left(  \varepsilon\left\Vert V;L^{2}\left(
\varpi_{+}\right)  \right\Vert +\varepsilon^{1/2}\left\Vert \varepsilon
^{1/2}R_{\varepsilon}^{-1}V;L^{2}\left(  \varpi_{+}\right)  \right\Vert
+\left\Vert V-\overline{V};L^{2}\left(  \varpi_{+}\right)  \right\Vert
\right)  \leq$\\
$\ \ \ \ \leq c\varepsilon^{3/2}\left\Vert V;\mathcal{H}_{\Omega}%
^{\varepsilon}\right\Vert =c\varepsilon^{3/2}.$%
\end{tabular}
\ \ \ \ \label{2.41}%
\end{equation}
For the first term (with $W_{i}$), we readily used the Schwarz inequality. For
the second term (with $1-X_{\varepsilon}$), we took into account that
$R\left(  x\right)  \leq c\varepsilon$ on supp$\left(  1-X_{\varepsilon
}\right)  $ and applied estimate (\ref{2.37}). For the third term (with
$V-\overline{V}$), we used estimate (\ref{2.36}).

Similar argument works for the last term $I_{3}^{\left(  p\right)  }$ in
(\ref{2.33}). Recalling boundary conditions in problem (\ref{2.7}) for the
asymptotic correctors $W_{i},$ we, indeed, obtain%
\[%
\begin{tabular}
[c]{l}%
$I_{3}^{\left(  p\right)  }=\left(  n^{\bullet}\cdot\nabla_{y}w^{\left(
p\right)  }+X_{\varepsilon}\underset{i=1}{\overset{2}{\sum}}\left(
\dfrac{\partial w^{\left(  p\right)  }}{\partial y_{i}}\partial_{n\left(
\xi\right)  }W_{i}+\varepsilon W_{i}n^{\bullet}\cdot\nabla_{y}\dfrac{\partial
w^{\left(  p\right)  }}{\partial y_{i}}\right)  +\right.  $\\
$\ \ \ \ \ \ \ \ \ \ \ \left.  +\varepsilon U_{\varepsilon}^{\left(  p\right)
}n^{\bullet}\cdot\nabla_{y}X_{\varepsilon},V\right)  _{\omega_{-}\left(
\varepsilon\right)  }=$\\
$\ \ \ \ \ =\left(  \left(  1-X_{\varepsilon}\right)  n^{\bullet}\cdot
\nabla_{y}w^{\left(  p\right)  }+\varepsilon U_{\varepsilon}^{\left(
p\right)  }n^{\bullet}\cdot\nabla_{y}X_{\varepsilon}+\varepsilon
X_{\varepsilon}\underset{i=1}{\overset{2}{\sum}}W_{i}n^{\bullet}\nabla
_{y}\dfrac{\partial w^{\left(  p\right)  }}{\partial y_{i}},V\right)
_{\omega_{-}\left(  \varepsilon\right)  }$%
\end{tabular}
\]
and, therefore,%
\begin{equation}%
\begin{tabular}
[c]{l}%
$\left\vert I_{3}^{\left(  p\right)  }-\varepsilon\left(  \underset
{i=1}{\overset{2}{\sum}}W_{i}n^{\bullet}\nabla_{y}\dfrac{\partial w^{\left(
p\right)  }}{\partial y_{i}},X_{\varepsilon}\overline{V}\right)  _{\omega
_{-}\left(  \varepsilon\right)  }\right\vert \leq$\\
$\ \ \ \ \ \ \ \ \leq c\left(  \left[  \text{meas}_{2}\left(  \text{supp}%
\left(  1-X_{\varepsilon}\right)  \right)  \right]  ^{1/2}\varepsilon
^{1/2}\left\Vert \varepsilon^{1/2}R_{\varepsilon}^{-1}V;L^{2}\left(
\omega_{-}\left(  \varepsilon\right)  \right)  \right\Vert +\right.  $\\
$\ \ \ \ \ \ \ \ \ \ \ \ \ \ \ \ \ \ \ \ \left.  +\varepsilon\left[
\underset{y\in\text{supp}X_{\varepsilon}}{\sup}R\left(  x\right)  ^{\rho
-1}\right]  \left\Vert V-\overline{V};H^{0}\left(  \omega_{-}\left(
\varepsilon\right)  \right)  \right\Vert \right)  \leq c\varepsilon
^{\min\left\{  1,\rho+1/2\right\}  }.$%
\end{tabular}
\ \ \label{2.42}%
\end{equation}
Here we applied inequalities (\ref{2.37}), (\ref{2.99}), (\ref{2.36}) while
taking into account the obvious relations%
\[
\text{meas}_{2}\left(  \text{supp}\left(  1-X_{\varepsilon}\right)  \right)
=O\left(  \varepsilon\right)  ,\underset{y\in\text{supp}X_{\varepsilon}}{\sup
}R\left(  x\right)  ^{\rho-1}=O\left(  \varepsilon^{\rho-1}\right)  ,\text{
\ }\rho\in\left(  0,1\right]  .
\]

It remains to examine the term%
\begin{equation}
I_{1}^{\left(  p\right)  }=-\left(  \Delta_{x}U_{\varepsilon}^{\left(
p\right)  },V\right)  _{\Omega_{\varepsilon}}, \label{2.43}%
\end{equation}
where, according to (\ref{2.13}),%

\begin{align*}
\Delta_{x}\mathbf{U}_{\varepsilon}^{\left(  p\right)  }  &  =\Delta
_{y}w^{\left(  p\right)  }\left(  y\right)  +\varepsilon U_{\varepsilon
}^{\left(  p\right)  }\left(  x\right)  \Delta_{y}X_{\varepsilon}%
+2\varepsilon\nabla_{y}U_{\varepsilon}^{\left(  p\right)  }\cdot\nabla
_{y}X_{\varepsilon}+\\
&  +X_{\varepsilon}\underset{i=1}{\overset{2}{\sum}}\left(  \varepsilon
^{-1}\Delta_{\xi}W_{i}\left(  \xi\right)  \dfrac{\partial w^{\left(  p\right)
}}{\partial y_{i}}\left(  y\right)  +2\nabla_{\eta}W_{i}\left(  \xi\right)
\cdot\nabla_{y}\dfrac{\partial w^{\left(  p\right)  }}{\partial y_{i}}\left(
y\right)  +\varepsilon W_{i}\left(  \xi\right)  \Delta_{y}\dfrac{\partial
w^{\left(  p\right)  }}{\partial y_{i}}\left(  y\right)  \right)  .
\end{align*}
Since $W_{i}$ is a harmonics, the first term in the sum vanishes. We now list
down estimates permitting to neglect some other terms:%

\begin{equation}%
\begin{tabular}
[c]{l}%
$\left\vert \left(  \left(  1-X_{\varepsilon}\right)  \Delta_{y}w^{\left(
p\right)  },V\right)  _{\Omega_{\varepsilon}}\right\vert \leq c%
{\displaystyle\int_{\Omega_{\varepsilon}\cap\text{supp}\left(
1-X_{\varepsilon}\right)  }}
R\left(  x\right)  ^{-1+\rho}\left\vert V\left(  x\right)  \right\vert dx\leq
$\\
$\ \ \ \ \ \leq c\left(
{\displaystyle\int_{-\varepsilon H}^{0}}
dz%
{\displaystyle\int_{0}^{c\varepsilon}}
R^{-2+2\rho}RdR\right)  ^{1/2}\varepsilon\left\Vert R_{\varepsilon}%
^{-1}V;L^{2}\left(  \Omega_{\varepsilon}\right)  \right\Vert \leq
c\varepsilon^{\rho+3/2},$\\
\\
$\varepsilon\left\vert \left(  U_{\varepsilon}^{\left(  p\right)  }\Delta
_{y}X_{\varepsilon}+2\nabla_{y}U_{\varepsilon}^{\left(  p\right)  }\cdot
\nabla_{y}X_{\varepsilon},V\right)  _{\Omega_{\varepsilon}}\right\vert \leq
c\varepsilon%
{\displaystyle\int_{\Omega_{\varepsilon}\cap\text{supp}\left\vert \nabla
_{y}X\varepsilon\right\vert }}
\left(  \varepsilon^{-1}+R\left(  x\right)  ^{-1+\rho}\right)  \left\vert
V\left(  x\right)  \right\vert dx\leq$\\
$\ \ \ \ \ \leq c\varepsilon\left(  \text{meas}_{3}\text{supp}\left\vert
\nabla_{y}X_{\varepsilon}\right\vert \right)  ^{1/2}\left\Vert R_{\varepsilon
}^{-1}V;L^{2}\left(  \Omega_{\varepsilon}\right)  \right\Vert \leq
c\varepsilon^{2},$\\
\\
$\varepsilon\left\vert \left(  X_{\varepsilon}W_{i}\Delta_{y}\dfrac{\partial
w^{\left(  p\right)  }}{\partial y_{i}},V\right)  _{\Omega_{\varepsilon}%
}\right\vert \leq c\varepsilon%
{\displaystyle\int_{\Omega_{\varepsilon}}}
X_{\varepsilon}\left(  y\right)  R\left(  x\right)  ^{-2+\rho}\left\vert
V\left(  x\right)  \right\vert dx\leq$\\
$\ \ \ \ \ \leq c\varepsilon\varepsilon^{-1+\rho}\left(  \text{meas}_{3}%
\Omega_{\varepsilon}\right)  ^{1/2}\left\Vert R_{\varepsilon}^{-1}%
V;L^{2}\left(  \Omega_{\varepsilon}\right)  \right\Vert \leq c\varepsilon
^{\rho+1/2},$\\
\\
$\left\vert \left(  X_{\varepsilon}\left(  \Delta_{y}w^{\left(  p\right)
}+2\underset{i=1}{\overset{2}{\sum}}\nabla_{\eta}W_{i}\cdot\nabla_{y}%
\dfrac{\partial w^{\left(  p\right)  }}{\partial y_{i}}\right)  ,V-\overline
{V}\right)  _{L^{2}\left(  \Omega_{\varepsilon}\right)  }\right\vert \leq$\\
$\ \ \ \ \ \ \ \ \ \ \ \ \ \ \ \ \ \ \ \ \ \ \ \ \ \ \leq c\varepsilon
^{-1+\rho}\left(  \text{meas}_{3}\Omega_{\varepsilon}\right)  ^{1/2}\left\Vert
V-\overline{V};L^{2}\left(  \Omega_{\varepsilon}\right)  \right\Vert \leq
c\varepsilon^{\rho+1/2}.$%
\end{tabular}
\ \ \label{2.44}%
\end{equation}
Here we have applied the same arguments as above.

Inequalities (\ref{2.44}) help to estimate all terms in (\ref{2.43}) with
exception of the first subtrahend on the left of (\ref{2.39}). Other two
subtrahends were exhibited in (\ref{2.42}) and (\ref{2.41}); hence, relation
(\ref{2.39}) is verified.

Similarly to the proof of Lemma \ref{Lemma2.4}, we now consider the integrals
over the cells $\Sigma_{\varepsilon}^{\nu}$ and their bases $\sigma
_{\varepsilon}^{\nu+}$ and $\sigma_{\varepsilon}^{\nu-},$ namely,%
\begin{equation}%
\begin{tabular}
[c]{l}%
$%
{\displaystyle\int_{\Sigma_{\varepsilon}^{\nu}}}
\Delta_{y}w^{\left(  p\right)  }\left(  y\right)  +\underset{i=1}{\overset
{2}{\sum}}\nabla_{\xi}W_{i}\left(  \dfrac{x}{\varepsilon}\right)  \cdot
\nabla_{\xi}\dfrac{\partial w^{\left(  p\right)  }}{\partial y_{i}}\left(
y\right)  X_{\varepsilon}\left(  y\right)  V\left(  y\right)  dx+$\\
$\ \ \ \ \ +\varepsilon%
{\displaystyle\int_{\sigma_{\varepsilon}^{\nu-}}}
W_{i}\left(  \dfrac{x}{\varepsilon}\right)  n^{\bullet}\left(  \dfrac
{x}{\varepsilon}\right)  \cdot\nabla_{y}\dfrac{\partial w^{\left(  p\right)
}}{\partial y_{i}}\left(  y\right)  X_{\varepsilon}\left(  y\right)  V\left(
y\right)  dy+\varepsilon\tau^{\left(  k\right)  }%
{\displaystyle\int_{\sigma_{\varepsilon}^{\nu+}}}
w^{\left(  p\right)  }yX_{i}\left(  y\right)  \overline{V}\left(  y\right)
ds_{x}.$%
\end{tabular}
\ \ \ \label{2.45}%
\end{equation}
Freezing here the argument $y$ at the mass center $y^{\nu}$ of the rectangle
$\sigma_{\varepsilon}^{\nu+},$ we perform integration in $\xi=\varepsilon
^{-1}x$ and, recalling calculation (\ref{2.9}), we obtain that expression
(\ref{2.45}) turns into the expression%
\[
\varepsilon^{3}\left(  \nabla_{y}\cdot b\nabla_{y}w^{\left(  p\right)
}\left(  y^{\nu}\right)  +\tau^{\left(  k\right)  }\left\vert \sigma
\right\vert w^{\left(  p\right)  }\left(  y^{\nu}\right)  \right)
X_{i}\left(  y^{\nu}\right)  \overline{V}\left(  y^{\nu}\right)
\]
which vanishes because $\left\{  \tau^{\left(  k\right)  },w^{\left(
p\right)  }\right\}  $ is an eigenpair of problem (\ref{2.10}), (\ref{2.11}).
The error of the freesing procedure does not exceed%
\begin{equation}
c\varepsilon%
{\displaystyle\int_{\sigma^{\nu+}}}
R\left(  y\right)  ^{-1+\rho}\left\vert X_{\varepsilon}\left(  y\right)
\overline{V}\left(  y\right)  -X_{\varepsilon}\left(  y^{\nu}\right)
\overline{V}\left(  y^{\nu}\right)  \right\vert dy \label{2.46}%
\end{equation}
while the weight factor $R^{-1+\rho}$ comes from inequality (\ref{2.99}) with
$k=2.$ Since%
\[%
\begin{tabular}
[c]{l}%
$X_{\varepsilon}\left(  y\right)  \overline{V}\left(  y\right)
-X_{\varepsilon}\left(  y^{\nu}\right)  \overline{V}\left(  y^{\nu}\right)
=$\\
$\ \ \ \ \ \ =%
{\displaystyle\int_{y_{1}^{\nu}}^{y_{1}}}
\dfrac{\partial\left(  X_{\varepsilon}\overline{V}\right)  }{\partial
\mathbf{y}_{1}}\left(  \mathbf{y}_{1},y_{2}^{\nu}\right)  d\mathbf{y}_{1}+%
{\displaystyle\int_{y_{2}^{\nu}}^{y_{2}}}
\dfrac{\partial\left(  X_{\varepsilon}\overline{V}\right)  }{\partial
\mathbf{y}_{2}}\left(  y_{1},\mathbf{y}_{2}\right)  d\mathbf{y}_{2},$%
\end{tabular}
\
\]
quantity (\ref{2.46}) may be bounded from above by
\[
c\varepsilon^{2}%
{\displaystyle\int_{\sigma_{\varepsilon}^{\nu+}}}
R^{-1+\rho}\left\vert \nabla_{y}\left(  X_{\varepsilon}\overline{V}\right)
\right\vert dy.
\]

Summing over $\nu_{i}\in\left[  -N_{i},N_{i}\right]  $ and applying relation
(\ref{2.35}), we get the following bound for the sum of all integrals
(\ref{2.46}):%
\begin{align*}
c\varepsilon^{2}%
{\displaystyle\int_{\omega}}
R^{-1+\rho}\left\vert \nabla_{y}\left(  X_{\varepsilon}\overline{V}\right)
\right\vert dy  &  \leq c\varepsilon^{1+\rho}%
{\displaystyle\int_{\omega}}
\left\vert \nabla_{y}\overline{V}\right\vert ^{2}+\left\vert \nabla
_{y}X_{\varepsilon}\right\vert ^{2}\left\vert \overline{V}\right\vert
^{2}dy\leq\\
&  \leq c\varepsilon^{1+\rho}%
{\displaystyle\int_{\omega}}
\left(  \left\vert \nabla_{y}\overline{V}\left(  y\right)  \right\vert
^{2}+R_{\varepsilon}^{-2}\left\vert \overline{V}\right\vert ^{2}\right)
dx\leq c\varepsilon^{\rho+1/2}\left\Vert V;\mathcal{H}_{\Omega}^{\varepsilon
}\right\Vert \leq c\varepsilon^{\rho+1/2}.
\end{align*}
Collecting the estimates obtained above, we conclude finally that
\begin{equation}
\delta=\left\Vert \mathcal{B}\mathbf{u}_{\varepsilon}^{\left(  p\right)
}-\mathbf{bu}_{\varepsilon}^{\left(  p\right)  };\mathcal{H}_{\Omega
}^{\varepsilon}\right\Vert \leq c_{k}\varepsilon^{\rho-1}. \label{2.47}%
\end{equation}

\subsection{The theorem on asymptotics of eigenvalues.}

Owing to (\ref{2.47}), Lemma \ref{Lemma 2.3} delivers an eigenvalue
$\beta_{\varepsilon}^{\left(  q\right)  }$ of the operator $\mathcal{B}%
_{\varepsilon}$ such that
\begin{equation}
\left\vert \beta_{\varepsilon}^{\left(  q\right)  }-\varepsilon^{-1}\left(
\tau^{\left(  k\right)  }+1\right)  ^{-1}\right\vert \leq c_{k}\varepsilon
^{\rho-1}. \label{2.50}%
\end{equation}
By (\ref{2.21}), we have $\beta_{\varepsilon}^{\left(  q\right)  }=\left(
\varepsilon+\alpha_{\varepsilon}^{\left(  q\right)  }\right)  ^{-1}$ and,
hence,%
\begin{equation}
\left\vert \alpha_{\varepsilon}^{\left(  q\right)  }-\varepsilon\tau^{\left(
k\right)  }\right\vert \leq c_{k}\varepsilon^{\rho-1}\varepsilon\left(
\tau^{\left(  k\right)  }+1\right)  \left(  \varepsilon+\alpha_{\varepsilon
}^{\left(  q\right)  }\right)  . \label{2.48}%
\end{equation}
This particularly gives
\[
\alpha_{\varepsilon}^{\left(  q\right)  }\leq\varepsilon\tau^{\left(
k\right)  }+c_{k}\varepsilon^{\rho+1}\left(  \tau^{\left(  k\right)
}+1\right)  +c_{k}\varepsilon^{\rho}\varepsilon\left(  \tau^{\left(  k\right)
}+1\right)  \alpha_{\varepsilon}^{\left(  q\right)  }.
\]
Thus, with a small $\varepsilon^{\left(  k\right)  }>0$ and $\varepsilon
\in\left(  0,\varepsilon^{\left(  k\right)  }\right]  $, we obtain%
\[
\alpha_{\varepsilon}^{\left(  q\right)  }\leq1-c_{k}\varepsilon^{\rho}\left(
\tau^{\left(  k\right)  }+1\right)  \varepsilon\left(  \tau^{\left(  k\right)
}+c_{k}\varepsilon^{\rho}\left(  \tau^{\left(  k\right)  }+1\right)  \right)
\leq c_{k}^{\alpha}\varepsilon
\]
and, according to (\ref{2.48}),%
\begin{equation}
\left\vert \alpha_{\varepsilon}^{\left(  q\right)  }-\varepsilon\tau^{\left(
k\right)  }\right\vert \leq c_{k}\varepsilon^{1+\rho}. \label{2.49}%
\end{equation}

\begin{theorem}
\label{Theorem 2.11} Let $\tau^{\left(  k\right)  }$ be an eigenvalue of
problem (\ref{2.10}), (\ref{2.11}) with multiplicity $\varkappa_{k},$ i.e.,
(\ref{2.24}) holds true. There exist $\varepsilon^{\left(  k\right)  }>0$ and
$C_{k}>0$ such that the eigenvalue sequence (\ref{1.88}) of problem
(\ref{1.12})-(\ref{1.14}) has at least $\varkappa_{k}$\ entries $\alpha
_{\varepsilon}^{\left(  Q\right)  },...,\alpha_{\varepsilon}^{\left(
Q+\varkappa_{k}-1\right)  }$ which satisfy inequality (\ref{2.49}).
\end{theorem}

\textbf{Proof.} It suffices to compute a bound for the number of eigenvalues
$\alpha_{\varepsilon}^{\left(  q\right)  }$ in (\ref{2.49}). Employing again
Lemma \ref{Lemma 2.3}, we set $\delta_{1}=Tc_{k}\varepsilon^{\rho-1}$ where
$c_{k}$ is taken from (\ref{2.47}). Then, for $p=k,...,k+\varkappa_{k}-1$, we
get coefficients $f_{\varepsilon^{q}}^{\left(  p\right)  }$ such that
\begin{equation}
\left\Vert \mathbf{u}_{\varepsilon}^{\left(  p\right)  }-\sum_{j=J\left(
\varepsilon\right)  }^{J\left(  \varepsilon\right)  +X\left(  \varepsilon
\right)  -1}f_{\varepsilon_{j}}^{\left(  p\right)  }u_{\varepsilon}^{\left(
j\right)  };\mathcal{H}_{\Omega}^{\varepsilon}\right\Vert \leq\frac{2\delta
}{\delta_{1}}\leq\frac{2}{T} \label{2.51}%
\end{equation}
where $\beta_{\varepsilon}^{\left(  J\left(  \varepsilon\right)  \right)
},...,\beta_{\varepsilon}^{\left(  J\left(  \varepsilon\right)  +X\left(
\varepsilon\right)  -1\right)  }$ is the list of all eigenvalues in
(\ref{2.17}) which meet the inequality%
\begin{equation}
\left\vert \beta_{\varepsilon}^{\left(  q\right)  }-\varepsilon^{-1}\left(
\tau^{\left(  k\right)  }+1\right)  ^{-1}\right\vert \leq Tc_{k}%
\varepsilon^{1+\rho}. \label{2.52}%
\end{equation}
Recall that the coefficient columns $f_{\varepsilon}^{\left(  p\right)
}=\left(  f_{\varepsilon J\left(  \varepsilon\right)  }^{\left(  p\right)
},...,f_{\varepsilon J\left(  \varepsilon\right)  +X\left(  \varepsilon
\right)  -1}^{\left(  p\right)  }\right)  ^{\intercal}\in\mathbb{R}^{X\left(
\varepsilon\right)  }$ are of unit length. Moreover, by (\ref{2.30}) and
(\ref{2.18}), we have%
\[%
\begin{tabular}
[c]{c}%
$\delta_{p,q}+O\left(  \varepsilon^{\min\left\{  \rho,1/2\right\}  }\right)
=\left\langle \mathbf{u}_{\varepsilon}^{\left(  p\right)  },\mathbf{u}%
_{\varepsilon}^{\left(  j\right)  }\right\rangle _{\varepsilon}=\left\langle
\sum\limits_{j=J\left(  \varepsilon\right)  }^{J\left(  \varepsilon\right)
+X\left(  \varepsilon\right)  -1}f_{\varepsilon_{j}}^{\left(  p\right)
}u_{\varepsilon}^{\left(  j\right)  },\sum\limits_{h=J\left(  \varepsilon
\right)  }^{J\left(  \varepsilon\right)  +X\left(  \varepsilon\right)
-1}f_{\varepsilon_{h}}^{\left(  q\right)  }u_{\varepsilon}^{\left(  h\right)
}\right\rangle +O\left(  T^{-1}\right)  =$\\
$=\left(  f_{\varepsilon}^{\left(  p\right)  }\right)  ^{\intercal
}f_{\varepsilon}^{\left(  q\right)  }+O\left(  T^{-1}\right)  .$%
\end{tabular}
\ \
\]
Thus, for small $\varepsilon$ and $T^{-1}$, the columns $f_{\varepsilon
}^{\left(  k\right)  },...,f_{\varepsilon}^{\left(  k+x_{k}-1\right)  }$ are
linear independent in $\mathbb{R}^{X\left(  \varepsilon\right)  }$ so that
$\varkappa_{k}\leq X\left(  \varepsilon\right)  .$ Since inequality
(\ref{2.52}) is just of the same kind as inequality (\ref{2.50}) which has
resulted in (\ref{2.49}), the proof of the theorem is completed.
\ \ $\blacksquare$

\section{Spectra of the problems\label{sect3}}

\subsection{Variational formulation of problems.\label{sect3.1}}

Let $\mathcal{H}$ denote the Sobolev space $H^{1}\left(  \Pi\left(
\varepsilon\right)  \right)  $ equipped with the specific norm%
\begin{equation}
\left\Vert \Phi;\mathcal{H}\right\Vert =\left(  \left\Vert \nabla_{x}%
\Phi;L^{2}\left(  \Pi\left(  \varepsilon\right)  \right)  \right\Vert
^{2}+\left\Vert \Phi;L^{2}\left(  \Lambda\left(  \varepsilon\right)  \right)
\right\Vert ^{2}\right)  ^{1/2} \label{3.1}%
\end{equation}
and the corresponding inner product (cf. (\ref{2.19})). We also introduce the
weighted Sobolev space $\mathcal{W}_{\theta}$ as the completion of
$C_{c}^{\infty}\left(  \overline{\Pi\left(  \varepsilon\right)  }\right)  $
(infinitely differentiable functions with compact supports) with respect to
the norm%
\begin{equation}
\left\Vert \Phi;\mathcal{W}_{\theta}\right\Vert =\left\Vert R_{\theta}%
\Phi;\mathcal{H}\right\Vert \label{3.2}%
\end{equation}
where $R_{\theta}=\exp\left(  \theta\left(  1+x_{1}^{2}\right)  ^{1/2}\right)
$ and $\theta\in\mathbb{R}.$ This space consists of all functions $\Phi\in
H_{loc}^{1}\left(  \overline{\Pi\left(  \varepsilon\right)  }\right)  $ with
the finite norm (\ref{3.2}). Clearly, $\mathcal{W}_{0}=\mathcal{H}.$ If
$\theta>0$, a function $\Phi\in\mathcal{W}_{\theta}$ decays exponentially as
$x_{1}\rightarrow\pm\infty$ but the space $\mathcal{W}_{\theta}$ with
$\theta<0$ includes functions with a certain exponential growth at infinity.

The standard formulation \cite{Lad} of the spectral problem (\ref{1.5}%
)-(\ref{1.7}) reads: to find $\lambda\in\mathbb{C}$ and $\Phi_{\varepsilon}%
\in\mathcal{H}\backslash\left\{  0\right\}  $ such that
\begin{equation}
\left(  \nabla_{x}\Phi_{\varepsilon},\nabla_{x}\Psi\right)  _{\Pi\left(
\varepsilon\right)  }=\lambda_{\varepsilon}\left(  \Phi_{\varepsilon}%
,\Psi\right)  _{\Lambda\left(  \varepsilon\right)  },\ \ \ \ \Psi
\in\mathcal{H}. \label{3.3}%
\end{equation}
For a fixed $\lambda,$ we also consider the integral identity
\begin{equation}
\left(  \nabla_{x}\Phi,\nabla_{x}\Psi\right)  _{\Pi\left(  \varepsilon\right)
}-\lambda\left(  \Phi,\Psi\right)  _{\Lambda\left(  \varepsilon\right)
}=F\left(  \Psi\right)  ,\ \ \ \ \Psi\in\mathcal{H}, \label{3.4}%
\end{equation}
serving for the inhomogeneous problem (\ref{1.5})-(\ref{1.7}) while
$F\in\mathcal{H}^{\ast}$ is a linear functional in the Hilbert space
$\mathcal{H}.$ A generalized solution of problem (\ref{1.5})-(\ref{1.7}) in
the weighted space $\mathcal{W}_{\theta}$ implies a function $\Phi
\in\mathcal{W}_{\theta}$ such that%
\begin{equation}
\left(  \nabla_{x}\Phi,\nabla_{x}\left(  R_{\theta}^{2}\Psi\right)  \right)
_{\Pi\left(  \varepsilon\right)  }-\lambda\left(  \Phi,R_{\theta}^{2}%
\Psi\right)  _{\Lambda\left(  \varepsilon\right)  }=F_{\theta}\left(
\Psi\right)  ,\ \ \ \Psi\in\mathcal{W}_{\theta}, \label{3.5}%
\end{equation}
where $F_{\theta}\in\mathcal{W}_{-\theta}^{\ast}.$ Formally, (\ref{3.5}) is
derived from (\ref{3.4}) by changing the test function $\Psi$ for the product
$R_{\theta}^{2}\Psi.$ Notice that the linear space $C_{c}^{\infty}\left(
\overline{\Pi\left(  \varepsilon\right)  }\right)  $ is dense in
$\mathcal{W}_{\theta}$ with any weight index $\theta.$

By the definition of the weighted norm (\ref{3.2}), $R_{\theta}^{2}\Psi
\in\mathcal{W}_{-\theta}$ in case $\Psi\in\mathcal{W}_{\theta}.$ Hence,
$\left(  \cdot,\cdot\right)  _{\Pi\left(  \varepsilon\right)  }$ and $\left(
\cdot,\cdot\right)  _{\Lambda\left(  \varepsilon\right)  }$ stand in
(\ref{3.5}) for extensions of the natural inner products in $L^{2}\left(
\Pi\left(  \varepsilon\right)  \right)  $ and $L^{2}\left(  \Lambda\left(
\varepsilon\right)  \right)  $ up to the duality between proper weighted
Lebesgue spaces.

Let $\theta>0$ and $F\in\mathcal{W}_{-\theta}^{\ast}\subset\mathcal{H}^{\ast}$
while $F_{\theta}\left(  \Psi\right)  =F\left(  R_{\theta}^{2}\Psi\right)  $
so that $F_{\theta}\in\mathcal{W}_{-\theta}^{\ast}$ as well. Then, if $\Phi
\in\mathcal{W}_{\theta}$ is a solution of problem (\ref{3.5}), $\Phi$ belongs
to $\mathcal{H}$ and is a solution of problem (\ref{3.4}) with an exponential
decay at infinity. Viceversa, in the case $\theta<0$ a solution $\Phi
\in\mathcal{H}$ of problem (\ref{3.4}), where $F\in\mathcal{H}^{\ast}%
\subset\mathcal{W}_{-\theta}^{\ast}$, becomes a solution of problem
(\ref{3.5}) in $\mathcal{W}_{\theta}.$

Under the symmetry assumption (\ref{1.8}), the same definition works for the
problem posed on the set (\ref{1.10}) with the artificial Dirichlet conditions
(\ref{1.11}). We use the notation $\mathcal{H}^{0}$ and \ $\mathcal{W}%
_{\theta}^{0}$ for the function space (\ref{1.9}) and the similar weighted
space of odd functions. Moreover, integral identities for this problem on
$\Pi_{\varepsilon}^{+}$ are refereed as the identities (\ref{3.3}),
(\ref{3.4}) and (\ref{3.5}) restricted onto the subspaces $\mathcal{H}^{0}$
and \ $\mathcal{W}_{\theta}^{0}$, respectively.

\subsection{The operator formulation of problems.\label{sect3.2}}

In the Hilbert space $\mathcal{H}$ with the inner product $\left\langle
\cdot,\cdot\right\rangle $, generated by norm (\ref{3.1}), we introduce the
operator $\mathcal{T}_{\varepsilon}$ by the formula
\begin{equation}
\left\langle \mathcal{T}_{\varepsilon}\Phi\text{ },\Psi\right\rangle =\left(
\Phi,\Psi\right)  _{\Lambda\left(  \varepsilon\right)  },\ \ \ \Phi,\text{
}\Psi\in\mathcal{H} \label{3.6}%
\end{equation}
(cf. formulae (\ref{2.17}) and (\ref{2.16}) in the domain $\Omega
_{\varepsilon}$). This operator is continuous with the unit norm, positive and
self-adjoint but not compact because the surface $\Lambda\left(
\varepsilon\right)  $ is unbounded. Thus its spectrum lies on the segment
$\left[  0,1\right]  $ of the real axis $\mathbb{R}\subset\mathbb{C}$ and its
essential spectrum does not reduce to the single point $\mu=0$ (see, e.g.,
\cite[Ch. 10]{BiSo}).

The restriction of $\mathcal{T}_{\varepsilon}$ on the subspace $\mathcal{H}%
^{0}$ is denoted by $\mathcal{T}_{\varepsilon}^{0}.$ Clearly, $\mathcal{T}%
_{\varepsilon}^{0}$ acts from $\mathcal{H}^{0}$ into $\mathcal{H}^{0}$. If
$\mu$ is an eigenvalue of the operator $\mathcal{T}_{\varepsilon}^{0}$ with
the eigenfunction $\Phi_{\varepsilon}^{0}\in\mathcal{H}^{0},$ analogously to
(\ref{2.21}),%
\begin{equation}
\lambda=\mu^{-1}-1 \label{3.7}%
\end{equation}
is an eigenvalue of problem (\ref{3.3}) restricted on $\mathcal{H}^{0},$ i.e.,
of the operator $\mathcal{L}_{\varepsilon}^{0}.$ Moreover, the odd extension
$\Phi_{\varepsilon}\in\mathcal{H}$ of the function $\Phi_{\varepsilon}^{0}$
over the plane $\left\{  x:x_{2}=0\right\}  $ becomes an eigenfunction of
problem (\ref{3.3}) (and problem (\ref{1.5})-(\ref{1.7}) on $\Pi\left(
\varepsilon\right)  $) corresponding to the same eigenvalue (\ref{3.7}). This
observation will be a tool to prove Theorem \ref{Theorem1.1}, the main result
of the paper.

Formula (\ref{3.7}) establishes a direct relation between the $\lambda
-$spectrum of the operator $\mathcal{L}_{\varepsilon}$ of problem (\ref{3.3})
and the $\mu-$spectrum of $\mathcal{T}_{\varepsilon}.$ Thus, we only examine
the spectra of the operators $\mathcal{T}_{\varepsilon}$ and $\mathcal{T}%
_{\varepsilon}^{0}$ in the sequel.

\subsection{Continuous spectra.\label{sect3.3}}

Clearly, the point $\mu=0$ is an eigenvalue of the operator $\mathcal{T}%
_{\varepsilon}$ with the infinite-dimensional eigenspace
\begin{equation}
\left\{  \Phi_{\varepsilon}\in\mathcal{H}:\text{\ }\Phi_{\varepsilon
}=0\text{\ on }\Lambda\left(  \varepsilon\right)  \right\}  \label{3.8}%
\end{equation}
A similar conclusion holds true for the operator $\mathcal{T}_{\varepsilon
}^{+}.$

\begin{lemma}
\label{Lemma3.1}The segment $\left(  0,1\right]  \subset\mathbb{R}%
\subset\mathbb{C}$ is filled with the continuous spectrum of the operator
$\mathcal{T}_{\varepsilon}.$
\end{lemma}

\textbf{Proof.} The assertion follows from general results \cite{Ko} (see also
\S 5.1 in \cite{NaPL}). For the reader convenience, we show here shortly how
to construct a singular Weyl sequence for any $\mu\in\left(  0,1\right]  $ so
that $\mu$ belongs to the essential spectrum of $\mathcal{T}_{\varepsilon}.$
Since the operator of problem (\ref{3.5}) regarded as the mapping
$\mathcal{W}_{\theta}\rightarrow\mathcal{W}_{-\theta}^{\ast}$ is Fredholm for
a sufficiently small negative $\theta$ (see \cite{Ko}, \cite[Theorem
5.1.4]{NaPL} and comments on the model problem (\ref{3.9}) below), the kernel
of this operator at $\theta=0,$ regarded as the mapping $\mathcal{H\rightarrow
H}^{\ast}$, is finite-dimensional. Thus, any point $\mu\in\left(  0,1\right]
$ of the essential spectrum lies in the continuous spectrum of $\mathcal{T}%
_{\varepsilon}.$

\begin{figure}
[ptb]
\begin{center}
\includegraphics[
height=1.6994in,
width=4.0352in
]%
{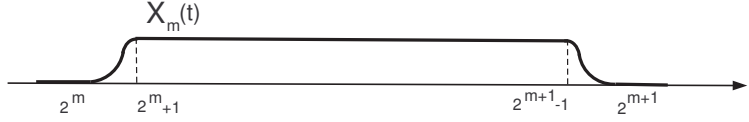}%
\caption{The plateau function.}%
\label{f6}%
\end{center}
\end{figure}

Let consider the model problem on the cross-section of the canal $\Pi,$ namely%
\begin{equation}%
\begin{tabular}
[c]{c}%
$-\Delta_{x^{\prime}}\varphi\left(  x^{\prime}\right)  +\eta^{2}\varphi\left(
x^{\prime}\right)  =0,$ \ \ \ \ $x^{\prime}\in\Gamma,$\\
$\partial_{x_{3}}\varphi\left(  x^{\prime}\right)  =\lambda\varphi\left(
x^{\prime}\right)  ,$ \ \ \ \ $x^{\prime}\in\gamma_{0},$ \ \ $\partial_{\eta
}\varphi\left(  x^{\prime}\right)  =0,$ \ \ \ \ $x^{\prime}\in\gamma.$%
\end{tabular}
\ \ \ \ \ \ \ \ \label{3.9}%
\end{equation}
Problem (\ref{3.9}) is obtained by the Fourier transform from the problem of
type (\ref{1.5})-(\ref{1.7}) in the cylindrical channel $\Pi=\mathbb{R}%
\times\Gamma$ while $\eta\in\mathbb{R}$ is the dual Fourier variable for
$x_{1}.$ Let $\mathcal{A}\left(  \lambda\right)  $ be an unbounded operator in
$L^{2}\left(  \Gamma\right)  $ associated (see \cite[Ch.10]{BiSo}) with the
bi-linear form%
\begin{equation}
Q\left(  \lambda,\varphi,\varphi\right)  =\left(  \nabla_{x^{\prime}}%
\varphi,\nabla_{x^{\prime}}\psi\right)  _{\Gamma}-\lambda\left(  \varphi
,\psi\right)  _{\gamma_{0}}. \label{3.10}%
\end{equation}
This operator is self-adjoint and bounded from below. Its domain belongs to
$H^{1}\left(  \Gamma\right)  .$ Since the embedding $H^{1}\left(
\Gamma\right)  \subset L^{2}\left(  \Gamma\right)  $ is compact, and
\begin{equation}%
\begin{tabular}
[c]{l}%
$Q\left(  \lambda_{1};\varphi,\varphi\right)  \geq Q\left(  \lambda
_{2};\varphi,\varphi\right)  ,$ \ \ $\lambda_{2}\geq\lambda_{1,}$%
\ \ $\varphi\in H^{1}\left(  \Omega\right)  ,$\\
$Q\left(  \lambda,1,1\right)  <0$ \ \ for $\lambda>0,$%
\end{tabular}
\ \ \ \ \ \ \ \ \label{3.11}%
\end{equation}
Theorems 10.1.2, 10.1.5, 10.2.4 in \cite{BiSo} ensure that the spectrum of
$\mathcal{A}\left(  \lambda\right)  $ is descrete and form the eigenvalue
sequence
\begin{equation}
\eta_{1}\left(  \lambda\right)  ^{2}<\eta_{2}\left(  \lambda\right)  ^{2}%
\leq\eta_{3}\left(  \lambda\right)  ^{2}\leq...\leq\eta_{k}\left(
\lambda\right)  ^{2}\leq...\rightarrow+\infty\label{3.12}%
\end{equation}
while $\mathbb{R}_{+}\ni\lambda\rightarrow\eta_{1}\left(  \lambda\right)
^{2}$ is a continuous, strictly monotone decreasing negative function. The
first eigenvalue $\eta_{1}\left(  \lambda\right)  ^{2}$ is simple due to the
maximum principle and $\eta_{1}\left(  \lambda\right)  =\pm i\left\vert
\eta_{1}\left(  \lambda\right)  \right\vert $ is imaginary. Let $\varphi
_{1}\left(  \lambda,x^{\prime}\right)  $ be the first eigenfunction of problem
(\ref{3.9}). We set%
\begin{equation}
\Phi^{\left(  m\right)  }\left(  x\right)  =a_{m}X_{m}\left(  \left(
2\pi\right)  ^{-1}\left\vert \eta_{1}\left(  \lambda\right)  \right\vert
x_{1}\right)  \sin\left(  \left\vert \eta_{1}\left(  \lambda\right)
\right\vert x_{1}\right)  \varphi_{1}\left(  \lambda,x^{\prime}\right)  ,
\label{3.13}%
\end{equation}
where $a_{m}$ is a normalization factor, $X_{m}$ is the plateau function in
Fig. \ref{f6},
\begin{equation}
X_{m}\left(  t\right)  =\chi\left(  t-2^{m}\right)  \chi\left(  2^{m+1}%
-t\right)  , \label{3.14}%
\end{equation}
and $\chi\in C^{\infty}\left(  \mathbb{R}\right)  $ is a cut-off function,
$\chi\left(  t\right)  =0$ for $t\leq0$ and $\chi\left(  t\right)  =1$ for
$t\geq1.$ The function $X_{m}$ is equal to one on the segment
\begin{equation}
\left[  2\pi\left\vert \eta_{1}\left(  \lambda\right)  \right\vert
^{-1}\left(  2^{m}+1\right)  ,2\pi\left\vert \eta_{1}\left(  \lambda\right)
\right\vert ^{-1}\left(  2^{m+1}-1\right)  \right]  \ni x_{1} \label{3.15}%
\end{equation}
and both functions (\ref{3.14}) and (\ref{3.13}) vanish for
\begin{equation}
x_{1}\notin\left[  2\pi\left\vert \eta_{1}\left(  \lambda\right)  \right\vert
^{-1}2^{m},2\pi\left\vert \eta_{1}\left(  \lambda\right)  \right\vert
^{-1}2^{m+1}\right]  . \label{3.16}%
\end{equation}
%

In the case $\lambda=0$ we simply set $\Phi^{\left(  m\right)  }\left(
x\right)  =a_{m}X_{m}\left(  x_{1}\right)  .$ We choose an integer $m$ such
that $\Theta\left(  \varepsilon\right)  \subset\left\{  x\in\Pi:x_{1}%
>2\pi\left\vert \eta_{1}\left(  \lambda\right)  \right\vert ^{-1}%
2^{m}\right\}  $ and obtain
\begin{equation}%
\begin{tabular}
[c]{l}%
$\left\langle \Phi^{\left(  m\right)  },\Phi^{\left(  m\right)  }\right\rangle
\geq a_{m}%
{\displaystyle\int_{2\pi\left\vert \eta_{1}\left(  \lambda\right)  \right\vert
^{-1}\left(  2^{m}+1\right)  }^{2\pi\left\vert \eta_{1}\left(  \lambda\right)
\right\vert ^{-1}\left(  2^{m+1}-1\right)  }}
\left(
{\displaystyle\int_{\Gamma}}
\left(  \left\vert \nabla_{x^{\prime}}\varphi_{1}\left(  \lambda;x^{\prime
}\right)  \right\vert ^{2}\left[  \sin\left(  2\pi\left\vert \eta_{1}\left(
\lambda\right)  \right\vert ^{-1}x_{1}\right)  \right]  ^{2}+\right.  \right.
$\\
$\ \ \ \left.  +\eta_{1}\left(  \lambda\right)  ^{2}\left\vert \varphi
_{1}\left(  \lambda,x^{\prime}\right)  \right\vert ^{2}\left[  \cos\left(
2\pi\left\vert \eta_{1}\left(  \lambda\right)  \right\vert ^{-1}x_{1}\right)
\right]  ^{2}dx^{\prime}+%
{\displaystyle\int_{\gamma_{0}}}
\left\vert \varphi\left(  \lambda;x^{\prime}\right)  \right\vert ^{2}%
dx_{2}\left[  \sin\left(  2\pi\left\vert \eta_{1}\left(  \lambda\right)
\right\vert ^{-1}x_{1}\right)  \right]  ^{2}dx_{1}\right)  =$\\
$\ \ =a_{m}^{2}\pi\left\vert \eta_{1}\left(  \lambda\right)  \right\vert
^{-1}\left(  2^{m+1}-2^{m}-2\right)  \left(  \left\Vert \nabla_{x^{\prime}%
}\varphi_{1};L^{2}\left(  \Gamma\right)  \right\Vert ^{2}+\left\vert \eta
_{1}\left(  \lambda\right)  \right\vert ^{2}\left\Vert \varphi_{1}%
;L^{2}\left(  \Gamma\right)  \right\Vert ^{2}+\left\Vert \varphi_{1}%
;L^{2}\left(  \gamma_{0}\right)  \right\Vert \right)  ^{2}.$%
\end{tabular}
\ \ \ \ \ \label{3.17}%
\end{equation}
We fix $a_{m}=O\left(  2^{-m/2}\right)  $ such that the last expression equals
$1.$ A similar calculation shows that $\left\langle \Phi^{\left(  m\right)
},\Phi^{\left(  m\right)  }\right\rangle $ is bounded from above uniformly in
$m$. The supports of the functions $\Phi^{\left(  m\right)  }$ and
$\Phi^{\left(  n\right)  }$ with $m\neq n$ are disjoint due to (\ref{3.16})
Hence, $\left\{  \Phi^{\left(  m\right)  }\right\}  $ converges to zero weakly
in $\mathcal{H}$ as $m\rightarrow+\infty$ and $\left\{  \Phi^{\left(
m\right)  }\right\}  $ implies a Weyl sequence provided
\begin{equation}
\left\Vert T_{\varepsilon}\Phi^{\left(  m\right)  }-\mu\Phi^{\left(  m\right)
};\mathcal{H}\right\Vert \rightarrow0\text{ \ as }m\rightarrow+\infty\text{
\ }. \label{3.18}%
\end{equation}
We have%
\begin{align*}
\left\Vert T_{\varepsilon}\Phi^{\left(  m\right)  }-\mu\Phi^{\left(  m\right)
};\mathcal{H}\right\Vert  &  =\sup\left\vert \left\langle T_{\varepsilon}%
\Phi_{m}^{\left(  m\right)  }-\mu\Phi^{\left(  m\right)  },\Psi\right\rangle
\right\vert =\\
&  =\mu\sup\left\vert \left(  \nabla_{x}\Phi^{\left(  m\right)  },\nabla
_{x}\Psi\right)  _{\Pi}-\lambda\left(  \Phi^{\left(  m\right)  },\Psi\right)
_{\Lambda}\right\vert =\\
&  =\mu\sup\left\vert -\left(  \Delta_{x}\Phi^{\left(  m\right)  }%
,\Psi\right)  _{\Pi}+\left(  \partial_{n}\Phi^{\left(  m\right)  }%
,\Psi\right)  _{\Lambda}\right\vert .
\end{align*}
Here the supremum is calculated over all $\Psi\in\mathcal{H}$ such that
$\left\Vert \Psi;\mathcal{H}\right\Vert =1.$ By the definition of $\eta
_{1}\left(  \lambda\right)  $ and $\varphi_{1}\left(  \lambda,\cdot\right)  $
as a solution of problem (\ref{3.19}), function (\ref{3.13}) satisfies the
boundary condition (\ref{1.7}) on $\Lambda$ and it is harmonic on a part of
the cylinder $\Pi$ determined by relation (\ref{3.15}). Thus, the expression
(\ref{3.19}) reduces to integral over the finite cylinders%
\[
\left\{  x\in\Pi:\left(  2\pi\right)  ^{-1}\left\vert \eta_{1}\left(
\lambda\right)  \right\vert x_{1}\in\left(  2^{m},2^{m}+1\right)  \right\}
\text{ \ and \ }\left\{  x\in\Pi:\left(  2\pi\right)  ^{-1}\left\vert \eta
_{1}\left(  \lambda\right)  \right\vert x_{1}\in\left(  2^{m+1}-1,2^{m}%
\right)  \right\}  .
\]
and, therefore, it does not exceed $ca_{m}\left\Vert \Psi;L^{2}\left(
\Pi\right)  \right\Vert \leq c2^{-m/2}.$ The proof is completed.
\ \ $\blacksquare$

The model problem on the cross section $\Gamma^{+}=\left\{  x^{\prime}%
\in\Gamma\text{ : }x_{2}>0\right\}  $ of cylinder (\ref{1.10})%
\begin{align}
-\Delta_{x^{\prime}}\varphi^{0}\left(  x^{\prime}\right)  +\eta^{2}\varphi
^{0}\left(  x^{\prime}\right)   &  =0,\ \ x^{\prime}\in\Gamma^{+}%
,\ \ \ \ \ \ \varphi^{0}\left(  x^{\prime}\right)  =0,\ \ x^{\prime}\in
\varpi^{0},\label{3.19}\\
\partial_{x_{3}}\varphi^{0}\left(  x^{\prime}\right)   &  =\lambda^{0}%
\varphi^{0}\left(  x^{\prime}\right)  ,\ \ x^{\prime}\in\gamma^{+}%
,\ \ \ \ \partial_{n}\varphi^{0}\left(  x^{\prime}\right)  =0,\ x^{\prime}%
\in\gamma_{0}^{+}\nonumber
\end{align}
corresponds to the operator $\mathcal{T}_{\varepsilon}^{0}$ of problem
(\ref{3.3}) restricted on $\mathcal{H}^{0}.$ Here $\varpi^{0}=\left\{
x^{\prime}\in\Gamma:x_{2}=0\right\}  $ and the curves $\gamma^{+},\gamma
_{0}^{+}$ compose the boundary of $\Gamma.$

First of all, we put $\eta=0$ and denote by $\lambda_{\Gamma}^{0}$ the first
eigenvalue of problem (\ref{3.19}) with the spectral Steklov boundary
condition. In view of the Dirichlet boundary condition, we have $\lambda
_{\Gamma}^{0}>0$ and
\begin{equation}
\mu_{\Gamma}^{0}=\left(  1+\lambda_{\Gamma}^{0}\right)  ^{-1}\in\left(
0,1\right)  , \label{3.20}%
\end{equation}
where $\varphi_{\Gamma}^{0}$ denotes the corresponding eigenfunction. Notice
that the trace inequality in $\Gamma^{+}$ reads:%
\begin{equation}
\left\Vert \varphi;L^{2}\left(  \gamma_{0}^{+}\right)  \right\Vert ^{2}%
\leq\left(  \lambda_{\Gamma}^{0}\right)  ^{-1}\left\Vert \nabla_{x^{\prime}%
}\varphi;L^{2}\left(  \Gamma^{+}\right)  \right\Vert ^{2},\text{\ }\varphi\in
H^{1}\left(  \Gamma^{+}\right)  ,\text{ \ }\varphi=0\text{\ on }\varpi^{0}.
\label{3.21}%
\end{equation}

We now examine the spectrum of the operator $\mathcal{T}_{\varepsilon}^{0}$ in
the space $\mathcal{H}^{0}$ which lies in the segment $\left[  0,1\right]  $
as well as the spectrum of $\mathcal{T}_{\varepsilon}.$

\begin{lemma}
\label{Lemma3.2}The segment $\left(  0,\mu_{\Gamma}^{0}\right]  $ is covered
with the continuous spectrum and the segment $\left(  \mu_{\Gamma}%
^{0},1\right]  $ contains descrete spectrum of the operator $\mathcal{T}%
_{\varepsilon}^{+}$ . The point $\mu=0$ is an eigenvalue with the
infinite-dimensional eigenspace $\left\{  \Phi_{\varepsilon}^{+}\in
\mathcal{H}^{0}:\Phi_{\varepsilon}^{+}=0\text{\ on }\Lambda\left(
\varepsilon\right)  \right\}  .$
\end{lemma}

\textbf{Proof.} To problem (\ref{3.19}) with the spectral parameter $\eta,$ we
associate the unbounded operator $\mathcal{A}^{0}\left(  \lambda\right)  $ in
the same way as in the proof of Lemma \ref{Lemma3.1}. If $\lambda
>\lambda_{\Gamma}^{0},$ the operator $\mathcal{A}^{0}\left(  \lambda\right)  $
meets the relation
\[
\left\langle \mathcal{A}^{0}\left(  \lambda\right)  \varphi_{\Gamma}%
^{0},\varphi_{\Gamma}^{0}\right\rangle =\left(  \lambda_{\Gamma}^{0}%
-\lambda\right)  \left(  \varphi_{\Gamma}^{0},\varphi_{\Gamma}^{0}\right)
_{\gamma_{0}^{+}}<0
\]
and, therefore, the first eigenvalue $\eta_{1}^{0}\left(  \lambda\right)
^{2}$ (cf. (\ref{3.12})) of problem (\ref{3.19}) with the fixed parameter
$\lambda>\lambda_{\Gamma}^{0}$ is negative. The same constructions as in
(\ref{3.13}) form the singular Weyl sequence. If $\lambda=\lambda_{\Gamma}%
^{0},$ we set $\Phi^{\left(  m\right)  }\left(  x\right)  =a_{m}X_{m}\left(
x_{1}\right)  \varphi_{\Gamma}^{0}\left(  x^{\prime}\right)  $ and again
conclude that point (\ref{3.20}) belongs to the continuous spectrum of
$\mathcal{T}_{\varepsilon}^{0}.$ It suffices to verify that $\left(
\mu_{\Gamma}^{0},1\right]  $ contains the descrete spectrum only. To this end,
we deal with the perturbed problem (\ref{3.4}) restricted on the subspace
$\mathcal{H}^{0},$ namely
\begin{equation}
\left(  \nabla_{x}\Phi^{+},\nabla_{x}\Psi\right)  _{\Pi^{+}\left(
\varepsilon\right)  }+\lambda M\left(  \Phi^{+},\Psi\right)  _{\Pi^{+}\left(
\varepsilon,L\right)  }-\lambda\left(  \Phi^{+},\Psi\right)  _{\Lambda
^{+}\left(  \varepsilon\right)  }=F^{+}\left(  \Psi\right)  ,\text{ }\Psi
\in\mathcal{H}^{0}, \label{3.22}%
\end{equation}
where $\Pi^{+}\left(  \varepsilon,L\right)  =\left\{  x\in\Pi^{+}\left(
\varepsilon\right)  :\left\vert x_{1}\right\vert <L\right\}  $ and $L$ is
chosen such that $\Theta\left(  \varepsilon\right)  \subset\left\{  x\in
\Pi\text{ }\left(  \varepsilon\right)  \text{: }\left\vert x_{1}\right\vert
<L\right\}  .$ Since the embedding $\mathcal{H}^{0}\subset L^{2}\left(
\Pi^{+}\left(  \varepsilon,L\right)  \right)  $ is compact, the difference of
the operator $\mathcal{T}_{\varepsilon}^{0}$ and the operator $\mathcal{T}%
_{\varepsilon,L}^{0},$ given by
\[
\left\langle \mathcal{T}_{\varepsilon,L}^{0},\Phi,\Psi\right\rangle =\left(
\Phi,\Psi\right)  _{\Lambda\left(  \varepsilon\right)  }-M\left(  \Phi
,\Psi\right)  _{\Pi\left(  \varepsilon,L\right)  },\text{ }\Phi,\Psi
\in\mathcal{H}^{0},
\]
is a compact operator. If we find $M$ such that the problem (\ref{3.22}) is
uniquely solvable for $\lambda\in\left[  0,\lambda_{\Gamma}^{0}\right)  $ and,
thus, the segment $\left(  \mu_{\Gamma}^{0},1\right]  $ is free of the
spectrum of $\mathcal{T}_{\varepsilon,L}^{0},$ then this segment contains only
the descrete spectrum of $\mathcal{T}_{\varepsilon}^{0}.$ To prove
additionally that a solution $\Phi^{+}\in\mathcal{H}^{0}$ of problem
(\ref{3.22}) decays exponentially at infinity, we transform the variational
problem (\ref{3.22}) into the following one which looks similar to
(\ref{3.5}):%
\begin{equation}
\left(  \nabla_{x}\Phi^{+},\nabla_{x}\left(  R_{\theta}^{2}\Psi\right)
\right)  _{\Pi^{+}\left(  \varepsilon\right)  }+\lambda M\left(  \Phi
^{+},R_{\theta}^{2}\Psi\right)  _{\Pi^{+}\left(  \varepsilon,L\right)
}-\lambda\left(  \Phi^{+},R_{\theta}^{2}\Psi\right)  _{\Lambda^{+}\left(
\varepsilon\right)  }=F^{+}\left(  R_{\theta}^{2}\Psi\right)  ,\text{ }\Psi
\in\mathcal{W}_{\vartheta}^{0}, \label{3.23}%
\end{equation}
where $F^{+}\in\left(  \mathcal{W}_{\theta}^{0}\right)  ^{\ast}.$ We put
$u=R_{\theta}\Phi^{+},$ $v=R_{\theta}\Psi$ and compute
\begin{align}
\left(  \nabla_{x}\Phi^{+},\nabla_{x}\left(  R_{\theta}^{2}\Psi\right)
\right)  _{\Pi^{+}\left(  \varepsilon\right)  }  &  =\left(  R_{\theta}%
\nabla_{x}\Phi^{+},\nabla_{x}v\right)  _{\Pi^{+}\left(  \varepsilon\right)
}+\left(  R_{\theta}\nabla_{x}\Phi^{+},vR_{\theta}^{-1}\nabla_{x}R_{\theta
}\right)  _{\Pi^{+}\left(  \varepsilon\right)  }=\label{3.24}\\
&  =\left(  \nabla_{x}u,\nabla_{x}v\right)  _{\Pi^{+}\left(  \varepsilon
\right)  }-\left(  uR_{\theta}^{-1}\nabla_{x}R_{\theta},\nabla_{x}v\right)
_{_{\Pi^{+}\left(  \varepsilon\right)  }}+\nonumber\\
&  +\left(  \nabla_{x}u,vR_{\theta}^{-1}\nabla_{x}R_{\theta}\right)  _{\Pi
^{+}\left(  \varepsilon\right)  }-\left(  uR_{\theta}^{-1}\nabla_{x}R_{\theta
},vR_{\theta}^{-1}\nabla_{x}R_{\theta}\right)  _{\Pi^{+}\left(  \varepsilon
\right)  }.\nonumber
\end{align}
Owing to the Lax-Milgram lemma, the inequality%
\begin{equation}
\left\Vert \nabla_{x}u;L^{2}\left(  \Pi\left(  \varepsilon\right)  \right)
\right\Vert ^{2}\leq cI\left(  u,u;\Pi\left(  \varepsilon\right)  \right)
\label{3.255}%
\end{equation}
for the left-hand side $I\left(  u,v;\Pi\left(  \varepsilon\right)  \right)  $
of (\ref{3.23}) provides the uniqueness and solvability of problem
(\ref{3.23}) together with the estimate of its solution%
\begin{equation}
\left\Vert \Phi^{+};\mathcal{W}_{\theta}^{0}\right\Vert \leq c\left\Vert
\nabla_{x}u;L^{2}\left(  \Pi\left(  \varepsilon\right)  \right)  \right\Vert
\leq c\left\Vert F;\left(  \mathcal{W}_{\theta}^{0}\right)  ^{\ast}\right\Vert
\label{3.25}%
\end{equation}
Let us prove (\ref{3.255}). First of all, we note that, for $\Psi=\Phi^{+},$
the second and third terms on the right of (\ref{3.24}) cancel each other.
Moreover,
\begin{equation}
\left\vert \nabla_{x}R_{\theta}\left(  x\right)  \right\vert \leq\theta
R_{\theta}\left(  x\right)  . \label{3.26}%
\end{equation}
Then we apply the trace inequality (\ref{3.21}) and the Friedrichs inequality%
\[
\left\Vert \varphi;L^{2}\left(  \Gamma^{+}\right)  \right\Vert ^{2}\leq
C\left\Vert \nabla_{x^{\prime}}\varphi;L^{2}\left(  \Gamma^{+}\right)
\right\Vert ^{2},\text{ \ \ }\varphi\in H^{1}\left(  \Gamma^{+}\right)
,\text{ \ \ }\varphi=0\text{\ on }\varpi^{0},
\]
both integrated over $x_{1}\in\mathbb{R}\backslash\left[  -L,L\right]  .$ As a
result, we obtain
\begin{equation}%
\begin{tabular}
[c]{l}%
$I\left(  u,u,\Pi^{+}\left(  \varepsilon\right)  \backslash\Pi^{+}\left(
\varepsilon,L\right)  \right)  \geq\left\Vert \nabla_{x}u;L^{2}\left(  \Pi
^{+}\left(  \varepsilon\right)  \backslash\Pi^{+}\left(  \varepsilon,L\right)
\right)  \right\Vert ^{2}-$\\
$\ \ \ \ \ \ \ \ \ \ \ \ -\theta^{2}\left\Vert u;L^{2}\left(  \Pi^{+}\left(
\varepsilon\right)  \backslash\Pi^{+}\left(  \varepsilon,L\right)  \right)
\right\Vert ^{2}-\lambda\left\Vert u;L^{2}\left(  \Lambda^{+}\backslash
\overline{\Pi^{+}\left(  \varepsilon,L\right)  }\right)  \right\Vert ^{2}\geq
$\\
$\ \ \ \ \ \ \geq\left(  1-\theta^{2}-\left(  \lambda_{\Gamma}^{0}\right)
^{-1}\lambda\right)  \left\Vert \nabla_{x}u;L^{2}\left(  \Pi^{+}\left(
\varepsilon\right)  \backslash\Pi^{+}\left(  \varepsilon,L\right)  \right)
\right\Vert ^{2}.$%
\end{tabular}
\ \ \ \ \ \ \ \ \ \label{3.27}%
\end{equation}
Finally, we use the similar three-dimensional inequalities on a finite part of
$\Pi\left(  \varepsilon\right)  $%
\begin{align*}
\left\Vert \Phi;L^{2}\left(  \Pi^{+}\left(  \varepsilon,L\right)  \right)
\right\Vert ^{2}  &  \leq\mathbf{c}\left(  \varepsilon,L\right)  \left\Vert
\nabla_{x}\Phi;L^{2}\left(  \Pi^{+}\left(  \varepsilon,L\right)  \right)
\right\Vert ^{2},\\
\left\Vert \Phi;L^{2}\left(  \Lambda^{+}\left(  \varepsilon\right)
\cap\overline{\Pi^{+}\left(  \varepsilon,L\right)  }\right)  \right\Vert ^{2}
&  \leq t\left\Vert \nabla_{x}\Phi;L^{2}\left(  \Pi^{+}\left(  \varepsilon
,L\right)  \right)  \right\Vert ^{2}+\mathbf{C}\left(  t,\varepsilon,L\right)
\left\Vert \Phi;L^{2}\left(  \Pi^{+}\left(  \varepsilon,L\right)  \right)
\right\Vert ^{2},
\end{align*}
and we derive that%
\begin{equation}%
\begin{tabular}
[c]{l}%
$I\left(  u,u;\Pi^{+}\left(  \varepsilon,L\right)  \right)  \geq\left\Vert
\nabla_{x}u;L^{2}\left(  \Pi^{+}\left(  \varepsilon,L\right)  \right)
\right\Vert ^{2}-v^{2}\left\Vert u;L^{2}\left(  \Pi^{+}\left(  \varepsilon
,L\right)  \right)  \right\Vert ^{2}-$\\
$\ \ \ \ \ \ \ \ \ \ -\lambda\left(  \left\Vert u;L^{2}\left(  \Lambda
^{+}\left(  \varepsilon\right)  \cap\overline{\Pi^{+}\left(  \varepsilon
,L\right)  }\right)  \right\Vert ^{2}-M\left\Vert u;L^{2}\left(  \Pi
^{+}\left(  \varepsilon,L\right)  \right)  \right\Vert ^{2}\right)  \geq$\\
$\ \ \geq\left(  1-\theta^{2}\mathbf{c}\left(  \varepsilon,L\right)
-t\lambda\right)  \left\Vert \nabla_{x}u;L^{2}\left(  \Pi^{+}\left(
\varepsilon,L\right)  \right)  \right\Vert ^{2}+\lambda\left(  M-\mathbf{C}%
\left(  t,\varepsilon,L\right)  \right)  \left\Vert u;L^{2}\left(  \Pi
^{+}\left(  \varepsilon,L\right)  \right)  \right\Vert ^{2}.$%
\end{tabular}
\ \ \ \ \ \ \ \ \label{3.28}%
\end{equation}
Set $M=\mathbf{C}\left(  t,\varepsilon,L\right)  $ to annul the last term in
(\ref{3.28}) and choose $\left\vert \theta\right\vert $ and $t>0$ sufficiently
small. Since $\lambda\in\left[  0,\lambda_{\Gamma}^{0}\right)  ,$ both the
factors on norms of $\nabla_{x}u$ on the right of (\ref{3.27}) and
(\ref{3.28}) stay positive. Hence, inequality (\ref{3.255}) and estimate
(\ref{3.25}) are valid. In terms of the operator $\mathcal{T}_{\varepsilon
,L}^{0}$ the latter with $\theta=0$ means that
\[
\left\Vert \left(  \left(  \mathcal{T}_{\varepsilon,L}^{0}-\mu\right)
^{-1}\Phi^{+};\mathcal{H}^{0}\right)  \right\Vert \leq c\left(  \mu\right)
\left\Vert \left(  \Psi;\mathcal{H}^{0}\right)  \right\Vert
\]
for $\mu\in\left(  \mu_{\Gamma}^{0},1\right]  .$ Thus, the operator
$\mathcal{T}_{\varepsilon,L}^{0}-\mu$ is an isomorphism. \ $\blacksquare$

\begin{corollary}
\label{Corollary 3.4}An eigenfunction $\Phi^{+}\in\mathcal{H}^{0}$ of the
operator $\mathcal{T}_{\varepsilon}^{0},$ corresponding to an eigenvalue
$\mu\in\left(  \mu_{\Gamma}^{0},1\right]  ,$ satisfies problems (\ref{3.22})
and (\ref{3.25}) with the functional%
\begin{equation}
\Psi\mapsto F^{+}\left(  \Psi\right)  =\left(  \mu^{-1}-1\right)  M\left(
\Phi,\Psi\right)  _{\Pi^{+}\left(  \varepsilon,L\right)  }. \label{3.29}%
\end{equation}
This functional is continuous on the weighted space $\mathcal{W}_{-\theta}%
^{0}$ with any $\theta\in\mathbb{R}$ because the integration set
$\overline{\Pi^{+}\left(  \varepsilon,L\right)  }$ in (\ref{3.29}) is compact.
Thus, $\Phi^{+}\in\mathcal{W}_{\theta}^{0}$ with a small $\theta>0$ so that
$\Phi^{+}$ decays exponentially at infinity. $\blacksquare$
\end{corollary}

\subsection{Discrete and point spectra.\label{sect3.4}}

The operator $-\mathcal{T}_{\varepsilon}^{0}$ is semi-bounded from below and,
by Lemma \ref{Lemma3.2}, it has the discrete spectrum in the segment $\left(
-1,-\mu_{\Gamma}^{0}\right]  .$ Let order the corresponding eigenvalues:%
\[
-\mu_{\varepsilon}^{\left(  1\right)  }\leq-\mu_{\varepsilon}^{\left(
2\right)  }\leq...\leq-\mu_{\varepsilon}^{\left(  \mathcal{N}\right)  }.
\]
We cannot exclude the case $\mathcal{N=}0$ yet and $\mathcal{N=}+\infty$ is
also possible.

The max-min principle (see \cite[Theorem 10.2.2]{BiSo}) applied for the
operator $-\mathcal{T}_{\varepsilon}^{0}$, reads:%
\begin{equation}
-\mu_{\varepsilon}^{\left(  k\right)  }=\underset{\mathcal{E}_{k}%
\subset\mathcal{H}^{0}}{\max}\underset{\Psi\in\mathcal{E}_{k}\backslash
\left\{  0\right\}  }{\inf}\frac{\left\langle -\mathcal{T}_{\varepsilon}%
^{0}\Psi,\Psi\right\rangle }{\left\langle \Psi,\Psi\right\rangle }.
\label{3.30}%
\end{equation}
Here $\mathcal{E}_{k}$ is any linear subspace of co-dimension $k-1$,
i.e.,$\dim\left(  \mathcal{H}^{0}\ominus\mathcal{E}_{k}\right)  =k-1$ and, in
particular, $\mathcal{E}_{1}=\mathcal{H}^{0}.$

Accepting the symmetry assumption (\ref{1.8}), we may consider the spectral
problem (\ref{1.12})-(\ref{1.14}) on the half $\Omega_{\varepsilon}%
^{+}=\left\{  x\in\Omega_{\varepsilon}:x_{2}>0\right\}  $ of the thin periodic
plate (\ref{1.4}), while prescribing the Dirichlet condition on the surface
$\left\{  x\in\Omega_{\varepsilon}:x_{2}=0\right\}  $. Let
\begin{equation}
0<\alpha_{\varepsilon}^{\left(  1\right)  +}<\alpha_{\varepsilon}^{\left(
2\right)  +}\leq...\leq\alpha_{\varepsilon}^{\left(  k\right)  +}%
\leq...\rightarrow+\infty\label{3.31}%
\end{equation}
be the ordered eigenvalue sequence of the formulated problem on $\Omega
_{\varepsilon}^{+}.$ The corresponding eigenfunctions $u_{\varepsilon
}^{\left(  k\right)  +}$ satisfy the relation%
\begin{equation}
\delta_{p,q}=\left(  \nabla_{x}u_{\varepsilon}^{\left(  p\right)  +}%
,\nabla_{x}u_{\varepsilon}^{\left(  q\right)  +}\right)  _{\Omega
_{\varepsilon}^{+}}=\alpha_{\varepsilon}^{\left(  p\right)  +}\left(
u_{\varepsilon}^{\left(  p\right)  +},u_{\varepsilon}^{\left(  q\right)
+}\right)  _{\omega_{+}^{+}} \label{3.32}%
\end{equation}
where $\omega_{+}^{+}=\left\{  x=\left(  y,z\right)  :y\in\omega,\text{ }%
y_{2}>0,\text{ }z=0\right\}  $ is the upper base of $\Omega_{\varepsilon}%
^{+}.$

Let also%
\begin{equation}
0<\tau^{\left(  1\right)  +}<\tau^{\left(  2\right)  +}\leq...\leq
\tau^{\left(  k\right)  +}\leq...\rightarrow+\infty\label{3.33}%
\end{equation}
be the eigenvalue sequence of the Dirichlet problem for the equation
(\ref{2.10}) on $\omega^{+}=\left\{  y\in\omega:y_{2}>0\right\}  .$ Theorem
\ref{Theorem 2.11} applied to the problems mentioned above, warrants the
inequality%
\[
\left\vert \alpha_{\varepsilon}^{\left(  q\right)  +}-\varepsilon\tau^{\left(
k\right)  +}\right\vert \leq c_{k}\varepsilon^{1+\rho}%
\]
for $\varepsilon\in\left(  0,\varepsilon_{k}\right]  $ and $\varepsilon
_{k}>0,$ $c_{k}>0$ depend on the eigenvalue number $k.$

We fix $N$ and put $\widetilde{\varepsilon}_{N}=\min\left\{  \varepsilon
_{1},...,\varepsilon_{N}\right\}  $ and $\widetilde{c}_{N}=\max\left\{
c_{1},...,c_{N}\right\}  .$ Now, for any $\varepsilon\in\left(  0,\widetilde
{\varepsilon}_{N}\right]  $ and $k=1,2,...,$ there exists $q=q\left(
k\right)  $ such that $q\left(  k_{1}\right)  \neq q\left(  k_{2}\right)  $
for $k_{1}\neq k_{2}$ and
\begin{equation}
\alpha_{\varepsilon}^{\left(  q\left(  k\right)  \right)  +}\leq
\varepsilon\tau^{\left(  k\right)  +}+\widetilde{c}_{N}\widetilde{\varepsilon
}_{N}^{\ \rho}\text{ \ }. \label{3.34}%
\end{equation}
Clearly, $q\left(  k\right)  \geq k$ and, therefore, $q\left(  k\right)  $ can
be changed for $k$ in (\ref{3.34}).

We extend the eigenfunctions $u_{\varepsilon}^{\left(  1\right)  +},$
...,$u_{\varepsilon}^{\left(  N\right)  +}$ by zero from $\Omega_{\varepsilon
}^{+}$ on $\Pi^{+}\left(  \varepsilon\right)  $ and keep the notation for
these extensions. If $k\leq N,$ any subspace $\mathcal{E}_{k}$ in (\ref{3.30})
contains a non-trivial linear combination%
\[
\Psi=a_{1}u_{\varepsilon}^{\left(  1\right)  +}+...+a_{k}u_{\varepsilon
}^{\left(  k\right)  +}\text{ \ }.
\]
Thus, the infimum in (\ref{3.30}) does not exceed%
\begin{align}
\dfrac{\left\langle -\mathcal{T}_{\varepsilon}^{0}\psi,\psi\right\rangle
}{\left\langle \psi,\psi\right\rangle }  &  =\dfrac{-\left\Vert \psi
;L^{2}\left(  \Lambda^{+}\left(  \varepsilon\right)  \right)  \right\Vert
^{2}}{\left\Vert \nabla_{x}\psi;L^{2}\left(  \Pi^{+}\left(  \varepsilon
\right)  \right)  \right\Vert ^{2}+\left\Vert \psi;L^{2}\left(  \Lambda
^{+}\left(  \varepsilon\right)  \right)  \right\Vert ^{2}}=\label{3.35}\\
&  =\dfrac{-\underset{j=1}{\overset{k}{\sum}}a_{j}^{2}\left\Vert
u_{\varepsilon}^{\left(  j\right)  +};L^{2}\left(  \omega_{+}^{+}\right)
\right\Vert ^{2}}{\underset{j=1}{\overset{k}{\sum}}a_{j}^{2}\left(
\alpha_{\varepsilon}^{\left(  j\right)  +}+1\right)  \left\Vert u_{\varepsilon
}^{\left(  j\right)  +};L^{2}\left(  \omega_{+}^{+}\right)  \right\Vert ^{2}%
}\leq-\dfrac{1}{1+\alpha_{\varepsilon}^{\left(  k\right)  +}}\leq-\dfrac
{1}{1+\varepsilon\left(  \tau^{\left(  k\right)  }+\widetilde{c}_{N}%
\widetilde{\varepsilon}_{N}^{\rho}\right)  }.\nonumber
\end{align}
Here we have used formulae (\ref{3.32}) and (\ref{3.34}). Hence, from the
max-min principle (\ref{3.30}) it follows that
\[
-\mu_{k}\leq-\left(  1+\varepsilon\left(  \tau^{\left(  k\right)  }%
+\widetilde{c}_{N}\widetilde{\varepsilon}_{N}^{\rho}\right)  \right)  ^{-1}%
\]
and, owing to (\ref{3.7}),%
\begin{equation}
\lambda_{\varepsilon}^{\left(  k\right)  }=\left(  1+\mu_{\varepsilon
}^{\left(  k\right)  }\right)  ^{-1}\leq\varepsilon\left(  \tau^{\left(
k\right)  }+\widetilde{c}_{N}\widetilde{\varepsilon}_{N}^{\rho}\right)  ^{-1}
\label{3.36}%
\end{equation}
If $d>0$ and $N$ are given, we choose $\varepsilon\left(  d,N\right)  >0$ such
that $\varepsilon\left(  d,N\right)  \leq\widetilde{\varepsilon}_{N}$ and,
with $k=1,...,N$ and $\varepsilon\in\left(  0,\varepsilon\left(  d,N\right)
\right)  ,$ the bound in (\ref{3.36}) does not exceed $d.$ Then Theorem 10.2.2
in \cite{BiSo} ensures the existence of, at least, $N$ eigenvalues
$-\mu_{\varepsilon}^{\left(  k\right)  }\in\left[  -1,-\left(  1+d\right)
^{-1}\right)  $ of the operator $-\mathcal{T}_{\varepsilon}^{0}$ which, as has
been explained, belong to the point spectrum of the operator $\mathcal{T}%
_{\varepsilon}.$ Corresponding numbers $\lambda_{\varepsilon}^{\left(
k\right)  }\in\left(  0,d\right)  $ are nothing but eigenvalues of problem
(\ref{1.5})-(\ref{1.7}) . Theorem \ref{Theorem1.1} is proved. $\blacksquare$

\section{Concluding remarks\label{sect4}}

\begin{figure}
[ptb]
\begin{center}
\includegraphics[
height=1.0577in,
width=3.2742in
]%
{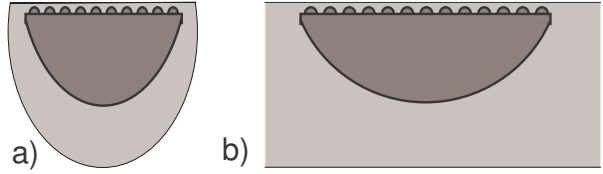}%
\caption{The transverse and longitudinal cross-sections of the body
$\Theta_{\mho}\left(  \varepsilon\right)  .$}%
\label{f7}%
\end{center}
\end{figure}

The main feature of the body $\Theta\left(  \varepsilon\right)  \subset
\Pi\left(  \varepsilon\right)  $ which provides the accumulation effect of the
trapped mode frequencies, is but the thin upper layer $\Omega_{\varepsilon}$
of water while the shape of the surface $\partial\Theta\left(  \varepsilon
\right)  \diagdown\partial\Omega_{\varepsilon}$ has no influence at all (cf.
\cite{NaMatSb} where $\omega_{-}\left(  \varepsilon\right)  =\partial
\Theta\left(  \varepsilon\right)  \diagdown\partial\Omega_{\varepsilon}$ is
flat). In this way, for the bodies $\Theta\left(  \varepsilon\right)  $ and
$\Theta_{\mho}\left(  \varepsilon\right)  $ with cross-section in Fig.
\ref{f3} and Fig. \ref{f7}, respectively, Theorem \ref{Theorem1.1} gives the
same bound $\varepsilon\left(  d,N\right)  $ for the small parameter
$\varepsilon$ in order to provide at least $N$ eigenvalues in the interval
$\left(  0,d\right)  $ of the continuous spectrum.%

If boundary of the periodic layer $\Sigma_{1}^{\infty}$ is Lipschitz only, the
convergence
\begin{equation}
\varepsilon^{-1}\alpha_{\varepsilon}^{\left(  q\right)  }\rightarrow
\tau^{\left(  q\right)  }\text{ \ for }\varepsilon\rightarrow0^{+} \label{5.1}%
\end{equation}
(cf. (\ref{2.49})) is valid. However, the homogenization technique to derive
(\ref{5.1}) differs from calculations performed in Section \ref{sect2} (cf.
\cite{BLP, SP} and others). Formula (\ref{5.1}) is sufficient to make the same
conclusion as in Theorem \ref{Theorem1.1}. In the estimate derived we
underline the convergence rate $O\left(  \varepsilon^{\rho}\right)  $ (cf.
(\ref{5.1}) and (\ref{2.49})) caused by singularities of $\nabla_{y}%
^{2}w\left(  y\right)  $ at the corner point of the rectangle $\omega$ (see
problem (\ref{2.10}), (\ref{2.11})).

In \cite[Ch. 7]{Nabook}, a method of inverse and direct reduction is developed
to describe an explicit dependence of constants $c_{k}$ in estimates of type
(\ref{2.49}) on the eigenvalue number $k$ and other attributes of the limit
spectrum (\ref{2.13}). This method requires rather intricate calculations and
we do not apply it here because the estimate (\ref{2.49}) is sufficient for
the main goal of the paper and the explicit dependence mentioned above does
not upgrade the result in Theorem \ref{Theorem1.1}.

\begin{figure}
[ptb]
\begin{center}
\includegraphics[
height=1.1096in,
width=2.8859in
]%
{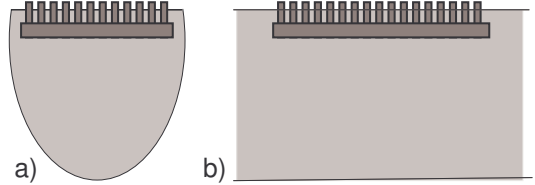}%
\caption{The transverse and longitudinal cross-section of the body
$\Theta_{\upuparrows}\left(  \varepsilon\right)  $ penetrating the surface.}%
\label{f8}%
\end{center}
\end{figure}

The shape of $\Theta_{\upuparrows}\left(  \varepsilon\right)  $ sketched on
Fig. \ref{f8}, where the rough surface $\partial\Theta_{\upuparrows}\left(
\varepsilon\right)  $ penetrates the water surface, is a possible
generalization. Although the plate $\Omega_{\varepsilon}$ becomes perforated,
it is very predictable that convergence (\ref{5.1}) and, thus, Theorem
\ref{Theorem1.1} remain valid though.%

\bigskip

\textbf{Acknowledgements.} This paper was prepared during the visit of S.A.
Nazarov to Department of Engineering of University of Benevento and to DIIMA
of University of Salerno and it was also supported by the grant RFFI-09-01-00759.

\end{document}